\documentclass[stsy,sglanonrev]{informs4}
\usepackage{eqndefns-left} % For checking the display equation width and equation environment definitions %
\RequirePackage{tgtermes}
\RequirePackage{newtxtext}
\RequirePackage{newtxmath}
\RequirePackage{bm}
\RequirePackage{endnotes}
\RequirePackage{amsmath}

\usepackage{mathtools}

\usepackage{graphicx} % Required for including images

\OneAndAHalfSpacedXII % Current default line spacing
\usepackage{algorithm}
\usepackage{algpseudocode}
\usepackage{tikz}

\usepackage{bbm}

\usepackage{todonotes}

\usepackage{natbib}
 \bibpunct[, ]{(}{)}{,}{a}{}{,}%
 \def\bibfont{\small}%

  \def\A{\cal A}%

\RequirePackage[english]{babel}
\RequirePackage[utf8]{inputenc}

\def\stochdomd{\preceq_{\textup{stoch}}}

\usepackage[font=it]{caption}

\EquationsNumberedThrough    % Default: (1), (2), ...
\TheoremsNumberedThrough     % Preferred (Theorem 1, Lemma 1, Theorem 2)
\ECRepeatTheorems  %  

\MANUSCRIPTNO{MOOR-0001-2024.00}

\begin{document}
%%%%%%%%%%%%%%%%

%\setcounter{page}{1}
%\DOI{10.1287/xxxx.0000.0000}%

\RUNAUTHOR{ Garcia, Bermolen, Jonckheere and Shneer}

% \RUNTITLE{Predictive Maintenance in Manufacturing}
\RUNTITLE{On Efficiency of Parallel and Restart Exploration}

% \TITLE{Optimal Resource Allocation in Humanitarian Logistics: A Stochastic Programming Approach}
%\TITLE{Parallelization and restart in stochastic simulations and applications to model-free settings}
\TITLE{Efficiency of Parallel and Restart Exploration Strategies in Model Free Stochastic Simulations}
%\TITLE {On Time Efficiency of Parallel and Restart Strategies for Stochastic Simulation}
%\TITLE{Probabilistic Insights for Efficient Exploration Strategies in Reinforcement Learning}

%\title{}
\ARTICLEAUTHORS{%
%\AUTHOR{John Doe,\textsuperscript{a} Jane Smith,\textsuperscript{b}}
%\AFF{\textsuperscript{a}Department of Industrial Engineering, University of XYZ, \EMAIL{john.doe@xyz.edu; \textsuperscript{b}Department of Computer Science, University of ABC, \EMAIL{jane.smith@abc.edu}} 

\AUTHOR{Ernesto Garcia}
\AFF{LAAS-CNRS, Toulouse, France. \EMAIL{egarciacig@laas.fr}}

\AUTHOR{Paola Bermolen}
\AFF{Facultad de Ingeniería, CICADA, Universidad de la República, Montevideo, Uruguay. \EMAIL{paola@fing.edu.uy}}

\AUTHOR{Matthieu Jonckheere}
\AFF{LAAS-CNRS, Toulouse, France. \EMAIL{matthieu.jonckheere@laas.fr}}

\AUTHOR{Seva Shneer}
\AFF{Heriot-Watt University, Edinburgh, UK. \EMAIL{V.Shneer@hw.ac.uk }}
% end of the block
}

\ABSTRACT{%
% Enter your abstract
We analyze the efficiency of parallelization and restart mechanisms for stochastic simulations in model-free settings, where the underlying system dynamics are unknown. Such settings are common in Reinforcement Learning (RL) and rare event estimation, where standard variance-reduction techniques like importance sampling are inapplicable. Focusing on the challenge of reaching rare states under a finite computational budget, we model exploration via random walks and Lévy processes. Based on rigorous probability analysis, our work reveals a phase transition in the success probability as a function of the number of parallel simulations: an optimal number $N^*$
  exists, balancing exploration diversity and time allocation per simulation. Beyond this threshold, performance degrades exponentially. Furthermore, we demonstrate that a restart strategy, which reallocates resources from stagnant trajectories to promising regions, can yield an exponential improvement in success probability. In the context of RL, these strategies can improve policy gradient methods by enabling more efficient state-space exploration, leading to more accurate policy gradient estimates. 
  }

\FUNDING{This research was supported by [SticAmsud LAGOON project, ANR EPLER]}

%Supplemental Material:
%Data Ethics & Reproducibility Note:

%\KEYWORDS{Stochastic programming, Decision support,Uncertainty, Disaster response, Optimization}

%\HISTORY{Received: Month DD, YYYY; Accepted: Month DD, YYYY; Published Online: Month DD, YYYY}

\maketitle

%%%%%%%%%%%%%%%%%%%%%%%%%%%%%%%%%%%%%%%%%%%%%%%%%%%%%%%%%%%%%%%%%%%%%%

% Text of your paper here

\section{Introduction}\label{sec:Intro}

%{\color{blue}
%- sparce or scarse\\
%- episodic context?
% }

The efficient simulation of stochastic systems is a cornerstone of many critical applications, ranging from rare-event estimation to gradient estimation in Reinforcement Learning. In many practical scenarios, the underlying dynamics of a system are complex, poorly characterized, or entirely unknown, placing us essentially in a model-free setting. In such contexts, strategies for exploring the state space where the stochastic dynamics live and gathering informative samples are paramount. This work investigates two powerful and general mechanisms for enhancing the efficiency of stochastic simulations: parallelization and restarting. We aim to provide a rigorous probabilistic analysis of these strategies, with applications focused on key model-free problems such as estimating small probabilities for stochastic processes with unknown dynamics and the challenge of efficient exploration in Reinforcement Learning (RL).

A canonical problem in probability and statistics is the estimation of the probability that a stochastic process reaches a rare, distant state within a finite time horizon. 
Standard Monte Carlo methods fail as the event of interest is rarely observed. A common technique for accelerating rare event simulation is importance sampling,
\cite{Asmussen_Glynn-Stochastic_simulation} which involves simulating the process under an alternative ``twisted" measure that makes the rare event more likely, and then correcting this change of measure via likelihood ratios. However, importance sampling is not a viable solution in the model-free settings as its application requires exact knowledge of the underlying process dynamics to construct the optimal change of measure. When the dynamics are unknown or too complex to model — as is often the case in large-scale RL environments or with intractable stochastic processes—this precise knowledge is unavailable, rendering importance sampling inapplicable. Moreover, an estimation scheme for the change of measure is usually of very large variance. It is precisely this limitation that motivates the use of blind or model-free strategies, such as parallelization and restarting, which do not require explicit knowledge of the system dynamics.

Parallelization—running multiple independent simulations concurrently—seems a natural remedy, as it increases the chances of observing the rare event. However, under a fixed total computational budget, a critical trade-off emerges: should one allocate the entire budget to a single, long simulation, or distribute it among many shorter, parallel runs? The answer is not trivial, as splitting the budget too finely may deprive each simulation of the necessary time to reach the rare state. We analyze this trade-off quantitatively, identifying a sharp phase transition in the success probability as a function of the number of parallel runs.

Beyond parallelization, restart strategies offer another avenue for improvement. Instead of letting simulations run in potentially unfruitful regions of the state space, a restart mechanism halts them and re-initializes them from more promising states, effectively redirecting computational effort. This restarting mechanism is particularly powerful in scenarios where certain regions are more likely to lead to the target event. A well-designed restart policy can prevent the waste of resources on trajectories with a low probability of success and serves as a practical, model-free alternative to importance sampling.

% in all, by leveraging parallel simulations and restart mechanisms, we can generate a more informative set of trajectories, leading to better gradient estimates and, consequently, more stable and efficient policy optimization in sparse-reward environments.

%However, efficient exploration appears as a critical challenge and often the primary bottleneck for the success of such algorithms (see \cite{Go_explore} and references therein). In many environments, the scarcity of meaningful signals throughout most of the simulation hinders the ability to conduct an effective statistical analysis of the reward function, which is crucial for policy updates. %When rewards are scarce, the design of appropriate control strategies that rely solely on the underlying process-often unknown-becomes essential for guiding the learning process. 

%{\color{blue} [VOLVER A INSISTIR SOBRE EL REGIMEN DE TIEMPO ACOTADO DE SIMULACIÓN PARA EVENTOS RAROS]
 Several recent works have addressed the efficient exploration challenge, in particular for RL. For example, the Go-Explore approach in \cite{Go_explore} introduces an innovative solution by distinctly separating the exploration and exploitation phases. This method allows the algorithm to maintain a repository of promising states, which can be revisited for further exploration, significantly enhancing the ability to discover and leverage rare or hard-to-reach states. 
%In contexts where a simulator or emulator is available, parallel simulation techniques also provide a powerful framework for enhancing pure exploration by running multiple instances of the environment simultaneously, each potentially following different strategies. This parallelism increases the chances of encountering rare states and enables a more exhaustive search of the state space. 
%Together with restarting mechanisms, these parallel approaches can further enhance the exploration process, significantly improving the overall effectiveness of RL algorithms in sparse reward environments.
Another example is found in \cite{FV_Jonckheere_etal} based on the paradigm of Fleming-Viot particle system, introduced by \cite{Burdzy_FV_original}. In this framework, a population of particles (which can represent different simulations in parallel) evolves over time, exploring the state space in parallel. When a particle falls into a less promising region, it is ``restarted" by being replaced with a copy of another, more promising particle. This ensures that resources are concentrated on exploring the most fruitful areas of the state space. 
%The Fleming-Viot strategy is particularly effective in scenarios where certain regions of the state space are more likely to yield valuable rewards, as it dynamically reallocates exploration effort towards these regions, thereby increasing the overall efficiency of the exploration process.

Beyond RL, stochastic restarting has been studied extensively across disciplines such as statistical physics, computer science, and network theory, primarily to minimize expected task-completion times modeled via first-passage times (see, e.g., \cite{Evans_stoch_reset, Luby_etal-Las_vegas, Avrachenkov_restart}). In contrast, our work focuses on the asymptotic behavior of the full distribution of passage times, studying large-deviation regimes and rare-event phenomena through asymptotic estimates, bounds, and threshold effects rather than expectations. While related analytic frameworks for Markov processes with restart—addressing ergodicity, quasi-stationary distributions, and connections to Fleming–Viot systems—are developed in \cite{Grigorescu_restart_Markov}, our results concern finer asymptotics for more restricted process classes and different probabilistic questions. Our work also differs from \cite{Monthus_restart_LDP}, which analyzes large deviations for the long-time asymptotics with stochastic restart to a fixed state. In contrast, we consider large deviations with respect to a spatial parameter modeling the complexity of exploring the state space.

Despite the promise of parallel simulation techniques and restart mechanisms, there is a lack of quantitative evidence of their potential benefits. Specifically, the extent to which these strategies improve exploration efficiency and yield better learning outcomes remains unclear. We aim to fill this gap by providing a detailed analysis of the impact of parallel simulations combined with restart mechanisms on the exploration process, modeled as the task of reaching a specific subset of the state space. To achieve this, we use a simplified model based on toy dynamics, specifically random walks and Lévy processes. These processes are chosen for their mathematical tractability and their ability to capture fundamental aspects of the problem, allowing explicit analytical results. Although these dynamics are simple, they provide crucial insights into the mechanisms that govern exploration, forming the basis for future work with more involved Markovian dynamics. 
By employing this framework, we examine two key aspects of exploration:

%By employing this framework, we focus on two key aspects of exploration:

\begin{itemize} \item {\bf Exploration complexity}, encapsulated by a parameter that quantifies the fluctuations necessary to achieve effective exploration and reach the target subset of the state space, 
\item {\bf  Exploration diversity}, represented by the number of parallel simulations employed, which increases the breadth of the state-space search and enhances the likelihood of encountering rare, high-reward states within the target subset. \end{itemize}

We focus our analysis on a specific regime where the exploration complexity scales linearly with the simulation budget. The toy dynamics framework allows us to examine these parameters within finite time budgets systematically and to quantify how parallel simulations improve efficiency by leveraging diversity, while restarting mechanisms enable particles to avoid stagnation by redirecting effort toward more promising regions of the state space, ultimately facilitating the achievement of the exploration goal.

%We focus on simple, well-defined Markovian dynamics to gain a clear understanding of how these strategies perform under different conditions. By systematically varying key parameters of the problem, such as the difficulty of reaching rare states and the degree of exploration diversity, we seek to identify the scenarios in which parallelism with restarting significantly enhances exploration and those in which, surprisingly, it may provide little or no benefit.

Our approach relies on the systematic analysis of these strategies under finite time budgets, using analytical techniques based on exponential martingales to quantify the effects of parallel simulations and restarting mechanisms. % These tools are particularly suited to studying rare event probabilities and light-tailed distributions, as detailed in \cite[Chapter VI]{Asmussen_Glynn-Stochastic_simulation}.   
By addressing the interplay between exploration diversity and the power of parallelization and restarting mechanisms, our work provides novel insights into efficient exploration strategies and lays the groundwork for further investigations into more complex stochastic dynamics.  

Through this investigation, our goal is not only to demonstrate the potential utility of these techniques but also to delineate their limitations. Understanding these nuances is crucial for developing more efficient stochastic simulations in model-free settings such as RL, particularly in environments characterized by sparse rewards. Ultimately, our findings aim to provide actionable insights into when and how to leverage parallel simulations and restarting strategies to optimize exploration, thereby advancing the field in both theory and practice.

%\subsection*{Related work} {\color{blue} separated or integrated section?}
%Though used in very different contexts, stochastic restarting has been studied in statistical physics. In this literature, it involves periodically restarting a process to minimize the \textit{expected time} to task completion. This technique can be analogously applied to RL, where restarting simulations can help escape local optima and encourage the discovery of new states.
% The restart rate may be considered deterministic or random, however, it was established in \cite{Luby_etal-Las_vegas} for discrete time and in \cite{Reuveni_stoch_restar_general} for continuous time, that the expected completion time is minimized for a deterministic rate. Let us also mention that the analytic computations of the optimal parameter is in general unfeasible. These results do not apply in our setting since our results require estimates on distributions rather than on expectations.

\subsection*{Problem Statement and Contributions}  

%{\color{blue}[CAMBIAR A STOCHASTIC SIMULATIONS]
In this work, we model exploration as a one-dimensional stochastic process, often referred to as a particle, aiming to reach a target subset of the state space. This target, represented as a barrier, captures a potentially informative state that is difficult to attain, making it essential to devise strategies that maximize the probability of success within a finite time budget and to quantify the resulting improvements.  In the following we summarize our main findings and key take-away messages. 

%\seva{Above, isn't the expected time to reach the barrier infinite? There is a positive probability to never reach it, right?}

%\paragraph{Main Findings and Contributions}  

\begin{enumerate}
    \item \textbf{Phase Transition in Reaching Probabilities.}  
    We identify a striking phase transition in the probability of reaching the target subset of the state space as a function of the number of parallel simulations $N$. In Theorems \ref{thm:main_theorem_num_particles} and \ref{thm:main_theorem_num_particles_Levy}, we rigorously show that when the time budget scales linearly with the barrier level, there is a sharp transition in performance. Our analysis is based on large deviations techniques and the transition happens at a certain threshold which we present explicitly. Below this threshold, concentrating the entire budget on a single particle is optimal. If the budget, when split into $N$ parts, remains above the threshold, then distributing it among $N$ particles significantly improves the probability of success. However, using too many particles—thereby reducing the time available to each beyond the threshold—results in exponentially diminished performance. 

    %\seva{I rephrased the above slightly but still not convinced it is as clear as it could/should be. Will return to it on my next pass.}

    \item \textbf{Optimal Number of Parallel Simulations.}
    In view of the previous point, a critical contribution is then the determination of an optimal number of parallel simulations, $ N^* $, that balances the trade-off between exploration diversity and the time allocated to each particle. This optimal $ N^* $ depends intricately on the time budget and the large deviations characteristics of the particles dynamics, but in terms of the phase transition mentioned above it can be characterized as the largest integer for which the split budget exceeds the threshold.  %Examples are provided to demonstrate that efficient estimation of \( N^* \) is feasible in practice.

    \item \textbf{Exponential Boost via restarting strategy.}  
%{\color{red}[ENFOCAR EN EL CASO DE   MEDIDA ``CUALQUIERA" PERO DEJAR QSD Y MENCIONAR FV  (que no es qsd en el transitorio) DE NUEVO EN ESTE CONTEXTO ]
    To overcome stagnation and improve the probability of reaching the target, we introduce a restarting mechanism inspired by mechanisms described, for instance, in \cite{FV_Jonckheere_etal}. We study an idealized situation in which this mechanism replaces particles with a low likelihood of success and restarts them at new initial positions sampled from a given measure.
    Our theoretical results (Theorem \ref{thm:restart_main_result_general_measure}) demonstrate that this restarting strategy can substantially enhance the probability of success. The improvement is quantified by a factor proportional to the time budget and the exponential moments of the measure associated with the restarting mechanism. A particularly interesting  restarting measure is a quasi-stationary distribution (QSD) over a set representing high-probability regions for successful trajectories. For these restarting measures, we can obtain sharper asymptotics using different methods that exploit the QSD properties.   We note that this restarting strategy can be well approximated in practice by selection mechanisms such as the Fleming-Viot or Aldous-Flanary-Palacios algorithms (see \cite{Fv_AFP_fraiman_etal} and the references therein). These schemes enable practical implementation and adequate sampling of initial positions, making the restarting mechanism computationally feasible in real-world scenarios.
 
   % \seva{I noticed that there was an idea to cite something in the paragraph above but the actual reference is missing.}
Finally, from a practical point of view, we present contributions in two directions for queueing systems. The first one corresponds to the problem of estimating small probabilities. We analyze the interplay between the probability of exploration and the efficiency of the associated estimation method. For the particular case of an M/M/1 queue, we show that the variance of the vanilla Monte Carlo estimator is essentially determined by the typical exploration time. 
The second, with broader impact, addresses efficient exploration in RL. In RL environments with sparse rewards, the core of the efficient exploration challenge lies in the inability of algorithms such as policy gradient methods to perform effective policy updates (see \cite{Go_explore} and references therein). These methods estimate the gradient of the expected return with respect to the policy parameters. This estimation relies crucially on exploring the state space; if the agent fails to reach rewarding states, the gradient estimates are biased and ineffective. Our work addresses this specific issue: we aim to improve gradient estimation (hence policy evaluation under a fixed policy) not by changing the policy, but by enhancing the exploration process used to evaluate it. The M/M/1/K queue considered in \cite{FV_Jonckheere_etal} is reanalized here to illustrate our results.

\end{enumerate}

%\paragraph{Methodological Approach.}  

\subsection*{Organization of the paper}

In Section \ref{sec:models}, we introduce the models for the exploration process, which represent the dynamics of the exploration phase in RL. Like the rest of the paper, this section is divided into two parts. The first part analyzes the effects of parallel exploration, and the second part analyzes our restarting mechanism. Section \ref{sec:Results} presents our main results for the two scenarios described above, after some necessary preliminaries. Main results are discussed in generality and illustrated through numerical simulations for concrete instances. In Section \ref{sec:proofs_parallel}, we prove Theorems \ref{thm:main_theorem_num_particles} and \ref{thm:main_theorem_num_particles_Levy} related to parallel exploration. In Section \ref{sec:proof_restart_general_measure}
 we prove Theorem \ref{thm:restart_main_result_general_measure} for exploration with a general restarting measure and 
in Section \ref{sec:proof_restart} we prove Theorem \ref{thm:main_thm_QSD_Levy} for a QSD restarting measure. 
Finally, appendices are devoted to practical applications of our work.

\section{Exploration strategies}
\label{sec:models}

{
We now formally introduce the models for exploration, using random walks and L\'evy processes as simplified toy models. We focus on the probability of reaching a distant state $x$ from an initial position of $0$, with the drift component encoding the difficulty of exploration. The drift represents the inherent challenge in progressing towards the target state, where a stronger drift away from the target reflects a more difficult exploration task.
}

{
We consider scenarios where the total time budget scales with $x$, allowing us to examine how exploration strategies perform as the target's reaching difficulty increases. Two specific strategies are explored: the first investigates parallel exploration under a fixed time budget, while the second examines the effects of restarting the process when it crosses a predefined threshold. Although these models are simplified, they provide critical insights into the fundamental principles that govern the effectiveness of exploration strategies in RL, particularly in challenging environments where reaching an unlikely target state is the goal.
}

\ 

\subsection{\bf Strategy I: parallel exploration.}
Let $\{Z_i(t)\}_{t\in I}, i \in \{1,\dots,N\}$, where $I= \mathbb N$ or $I= \mathbb R^+$, be a family of
one-dimensional random walks or L\'evy processes.
Define further
\begin{align} 
\tau_i(x) &\coloneq \inf \{ t\in I:\ \ Z_i(t) \ge x\},  \\
\tau^{(N)}(x)& := \min \{ \tau_1(x), \tau_2(x), \dots, \tau_N(x) \},  
\label{def:tau}
\end{align}
the time it takes particle $i$ to reach state $x$, and the first time one of $N$ particles reaches $x$, respectively.

First, the case of $N$ parallel independent simulations is analyzed.
We introduce as performance criterion for exploration the probability that starting from $0$, one of the simulations reaches the rare state $x$ within a given time budget $B(x)/N$, i.e.
$$  \PP\Big( \tau^{(N)}(x) \leq B(x)/N \Big).$$

%\seva{Is it $B(x)/N$ above? here is the criterion and then it makes sense for budgets scaling correctly but maybe is confusing...}

We aim at characterizing this probability (as a function of $N$) in the regime where the total simulation budget, $B(x)$, satisfies $ x = C\cdot B(x) + o \left(\sqrt{B(x)}\right)$, for $x$ large and $C$ some positive constant. This is the focus of Theorem \ref{thm:main_theorem_num_particles} for random walks and Theorem \ref{thm:main_theorem_num_particles_Levy} for L\'evy processes.

\subsection{\bf Strategy II: exploration with restart.}

Our second model incorporates a \textit{restart mechanism} into the dynamics of a random walk or a L\'evy process. The ingredients of the model are: a positive number $x>0$ (the rare state we want the process to reach), a probability measure $\nu_x$ supported on $(0,x)$, and a discrete random walk or a L\'evy process $Z = \{Z(t)\}_{t\in\RR^+}$. The trajectories of the process under study can be viewed as a concatenation in time of independent and identically distributed cycles, each consisting of sampling a state $y_0$ from $\nu_x$ and simulating the trajectory of $y_0+Z(t)$ until the first time it exits the interval $(0,x)$.

Our goal is to analyse the first passage time over $x$ of the restarted process,
% We intend to understand the time it takes for the process to exit through $[x,+\infty)$, namely, the first passage time over $x$ of the restarted process, and 
and compare it with the corresponding one for the process without restart. 
The conditions required on the family of restart measures 
to improve algorithmic performance are quite general. Essentially, it suffices that each $\nu_x$ be stochastically dominated by a non-degenerate measure possessing a finite second moment. 
%The conditions on the restarting family measures $\nu_x$ to boost the performance of the algorithm are pretty general: essentially they should be stochastically dominated by a non-degenerate measure with finite second moment.
We also examine in detail the case where restart occurs according to the unique quasi-stationary measure of the process 
%We also analyze the particular case of the unique quasi-stationary measure for the process 
$Z$ absorbed at $\RR^+\setminus (0,x)$  (see Section \ref{subsec:Model_restart_qsd} for precise definitions). Our motivation for this particular mechanism (restarting from quasi-stationary distributions) is that the resulting process approximates complex particle dynamics built upon a given Markov process, a key example, as noted earlier, being the class of Fleming–Viot particle systems. 
%The example we have in mind is the already mentioned Fleming-Viot particle systems. 
%, which consists, given a Markov process with state space $\cS$, a number of particles $N\in\NN$ and a subset $A\subset \cS$, of iterating the following procedure: $N$ independent copies of the process are simulated on $\cS\setminus A$ until the time any of them attempts to enter $A$. The trajectory of that particle is halted and immediately restarted at the current position of one of the others, chosen at random.  
%In \cite{FV_Jonckheere_etal} the authors explore the use of FV dynamics for learning purposes, both for exploration and estimation in different contexts. Our study is intended as an initial step in the development of a theoretical framework around it. Although simple to describe, a rigorous mathematical treatment of the FV dynamics is a challenging task. 
Our model is indeed an approximation of the Fleming-Viot dynamics since the stationary distribution of the $N-$ particle system approximates the quasi-stationary distribution associated with the underlying Markov process absorbed upon reaching $A$, as $N$ goes to infinity (see 
%for example \cite{Asselah_Ferrari_Groisman_QSd_finite_spaces, ferrari_maric, Villemonais_FV_to_QSD_diffusions} and 
Appendix II for more details).

\section{Main results}\label{sec:Results}

After introducing due notation, we present our results for each strategy.

\subsection{Parallel exploration driven by a random walk}\label{subsec:Model_parallel}

Let $Z(t) = X_1+\dots+X_t$, $t\in\NN$ be a random walk with independent increments all distributed as a random variable $X$. We will always assume that $X$ is not deterministic, meaning that it is not concentrated on a single value, and that it may take both positive or negative values with positive probability. The moment-generating function (mgf) of $X$, $\varphi(\lambda):= \EE[\exp(\lambda X)]$, is strictly convex on the interior points of the set $\Lambda$ ($\text{int}(\Lambda)$) where
$
\Lambda := \{\lambda\in \RR:\ \varphi(\lambda) <+\infty\}.
$
%Let also $\psi(\lambda)\coloneqq \log \varphi(\lambda)$. 
The following condition is usually referred to as the \textit{right Cram\'er condition:}
\begin{equation}\label{right_cramer_cond}
  \Lambda_+ := (0,+\infty)\cap \text{int}(\Lambda)\ \textit{is non empty,}
\end{equation}
We will work under further additional assumptions:
\begin{equation}\label{diff_mgf}
        \textit{ $\varphi$ is at least twice differentiable in } \Lambda_+,
\end{equation} and
\begin{equation}\label{positive_cramer_exp}
   \textit{There exists } \lambda^*  \in \Lambda_+:\ \psi(\lambda^*) = 0 \quad \text{ where } \quad \psi(\lambda)\coloneqq \log \varphi(\lambda)
\end{equation}
%If $\EE[|X|^k]<+\infty$ then $\varphi^{(k)}(0)^+ = \EE[X^k]$, $k\geq 1$, where $\varphi^{(k)}(\cdot)^+$ denotes the right derivative of order $k$. 

Note that $\EE[X] = \psi'(0)$ and assumptions \eqref{right_cramer_cond}, \eqref{diff_mgf} and \eqref{positive_cramer_exp} imply that $\EE[X] < 0$ and $\psi'(\lambda^*) > 0$, since $\psi(\lambda^*) = \psi(0)$ and $\varphi$ is strictly convex. 
% We now introduce further notation required for the statement of our main result. 
Let \( \tau(x) \) denote the first time the random walk reaches level \( x \), defined as:
 \begin{equation}\label{def:passage_time}
 \tau(x) := \inf\{ t \geq 0 : Z(t) \geq x \}.    
 \end{equation} Furthermore, let \( \tau^{(N)}(x) \) represent the minimum of \( N \) independent random times, each distributed as \( \tau(x) \), i.e.,
 \begin{equation}\label{def:min_of_passage_times}
 \tau^{(N)}(x) := \min \{ \tau_1(x), \tau_2(x), \dots, \tau_N(x) \},    
 \end{equation} which corresponds to the first time any of \( N \) independent random walks reaches level \( x \). 

Our first main theorem states that the interesting behaviour of $\tau(x)$ and $\tau^{(N)}(x)$ occurs while comparing them to functions of the following classes: for $\lambda\in\Lambda_+$ let \begin{equation}\label{def:asymp_linear_funcs}
    L(\lambda) \coloneqq \left\{f:\RR\to\RR\ \text{ such that } x = \psi'(\lambda)\cdot f(x) + \textbf{o}\left(\sqrt{f(x)}\right)\ \text{as } x\to+\infty\right\}.
\end{equation} This set consists of functions \( f(x) \) that asymptotically behave like \( x/\psi'(\lambda) \) with an error term of order smaller than \( \sqrt{f(x)} \).
The result below is stated for the class of functions defined above but the reader can think of a specific example being $B(x) = x/\psi'(\lambda)$ for any $\lambda\in\Lambda_+$. 

We now state our main result regarding parallel exploration for random walks:

\begin{theorem}\label{thm:main_theorem_num_particles}
    Let $Z$ denote a random walk whose increments satisfy assumptions \eqref{right_cramer_cond}, \eqref{diff_mgf} and \eqref{positive_cramer_exp}, and let $B(\cdot)$ belong to $L(\lambda)$ for some $\lambda\in\Lambda_+$. Then for $N\geq 2$:
    \begin{equation}\label{eq:threshold_number_particles}
        \lim_{x\to+\infty}\frac{\PP\left(\tau^{(N)}(x)\leq \frac{B(x)}{N}\right)}{\PP\left(\tau(x)\leq B(x)\right)} =  \begin{cases}
            N,\ \text{ if } N \psi'(\lambda) < \psi'(\lambda^*);\\
            0,\ \text{ if } N \psi'(\lambda) > \psi'(\lambda^*).
            %\\
            %0,\ \text{ if } \mu(\lambda^*)<\mu(\lambda)%<N\mu(\lambda).
        \end{cases}
    \end{equation}
\end{theorem}

\begin{corollary}[Optimal number of particles] \label{cor:optimal_num_part}
For a given state $x$ and a budget $B(x) \in L(\lambda)$. If $\psi'(\lambda) < \psi'(\lambda^*)$, the asymptotically optimal number of particles is:
\begin{equation}
    N^* = \max \; \left\{ N\geq 1 :\ N\psi'(\lambda) < \psi'(\lambda^*)\right\} =  \left\lceil \frac{\psi'(\lambda^*)}{\psi'(\lambda)} \right\rceil - 1,
\end{equation} where $\lceil\cdot\rceil$ is the ceiling function. Otherwise $N^*=1$.
\end{corollary}
In the regime in which the theorem is stated, probabilities decay exponentially as $x$ goes to infinity (see Proposition \ref{thm:estimates_pasasge_times} in Section \ref{sec:proofs_parallel}). The result shows that using too many particles (thus giving too little time to each) entails an exponentially worse performance, meaning that $\PP\left(\tau^{(N)}(x)\leq \frac{B(x)}{N}\right)$ decays with a faster exponential rate than $\PP\left(\tau(x)\leq B(x)\right)$. On the other hand using too few particles comes only at a linear cost in performance. 

%{\color{purple} 
The main intuition behind the results stated above is that within the (rare) event $\{\tau(x) < +\infty\}$, the probability is concentrated on those trajectories with drift $\psi'(\lambda^*)$, for which the time needed to reach $x$ is $\tau(x) \approx x/\psi'(\lambda^*)$. This idea is formalized using exponential martingales and is at the core of Proposition \ref{thm:estimates_pasasge_times}. Our results show that allowing a budget smaller than $x/\psi'(\lambda^*)$ substantially reduces the probability of success, while allowing too much budget results in no significant improvement. In the latter case, by splitting it into several independent particles, each one provides an ``extra chance" of success.
%}

%\seva{This is perhaps something to discuss in a meeting but at the moment I'm not sure the above in purple explains things clearly. One could try to expand it, starting from the very last observation on the expectation of the time to reach $x$, or move it elsewhere. I'm not sure what to do at the moment.}

Even though in general $\lambda^*$ and (thus) the threshold $\psi'(\lambda^*)^{-1}$ depend in a non-trivial way on the distribution of the increments, in some cases it is possible to find simple expressions for these parameters. We present some of the important examples here.

\begin{example}\textbf{Normal increments:}

\label{example:gaussian_increments_rw}
    If $X$ has a normal distribution with mean $-\mu<0$ and variance $\sigma^2=1$, then: \begin{equation*}
        \psi(\lambda) = -\lambda \mu + \frac{1}{2}\lambda^2,
    \end{equation*} hence $\lambda^* = 2\mu$ and $\psi'(\lambda^*) = \mu$ and the threshold can be reliably estimated in practice if needed. Moreover, the function $\psi':\Lambda_+ \to [0,+\infty)$ is surjective, so for every $C>0$ there exists $\lambda\in\Lambda_+$ such that $\psi'(\lambda) = 1/C$. Theorem \ref{thm:main_theorem_num_particles} in this case states that if, in order to reach a high state $x>0$, we are allowed a budget $B(x) = Cx$ for some $C>0$, then the (asymptotically) optimal number of independent particles to use is $N^* = \lceil C \mu\rceil - 1$ if $C\mu>1$ and $N^* = 1$ if not.
    \rqed
\end{example}

\begin{example}\label{example:BD_chains}\textbf{Birth-and-death chains:}
    If $X$ only takes values $1$ and $-1$ with probabilities $p$ and $1-p$ respectively with $p<1/2$, then $\lambda^* = \log(1-p) - \log p$ and $\psi'(\lambda^*) = 1-2p = -\EE X$. Since $X$ is bounded, the image of $\psi'$ contains the positive real numbers as in the previous example.  Then if we are allowed a budget $B(x) = Cx$, the number of particles maximizing our chances is $N^* = \lceil C (1-2p)\rceil - 1$ if $C(1-2p)>1$ and $N^* = 1$ if not.
    \rqed
\end{example}

%\seva{In these examples, can we say what the theorem boils down to? I.e., can we present optimal $N$ in these cases?}

\subsubsection{Simulations for parallel random walks:}

\begin{figure}
    \centering
    \includegraphics[width = 0.9\linewidth]{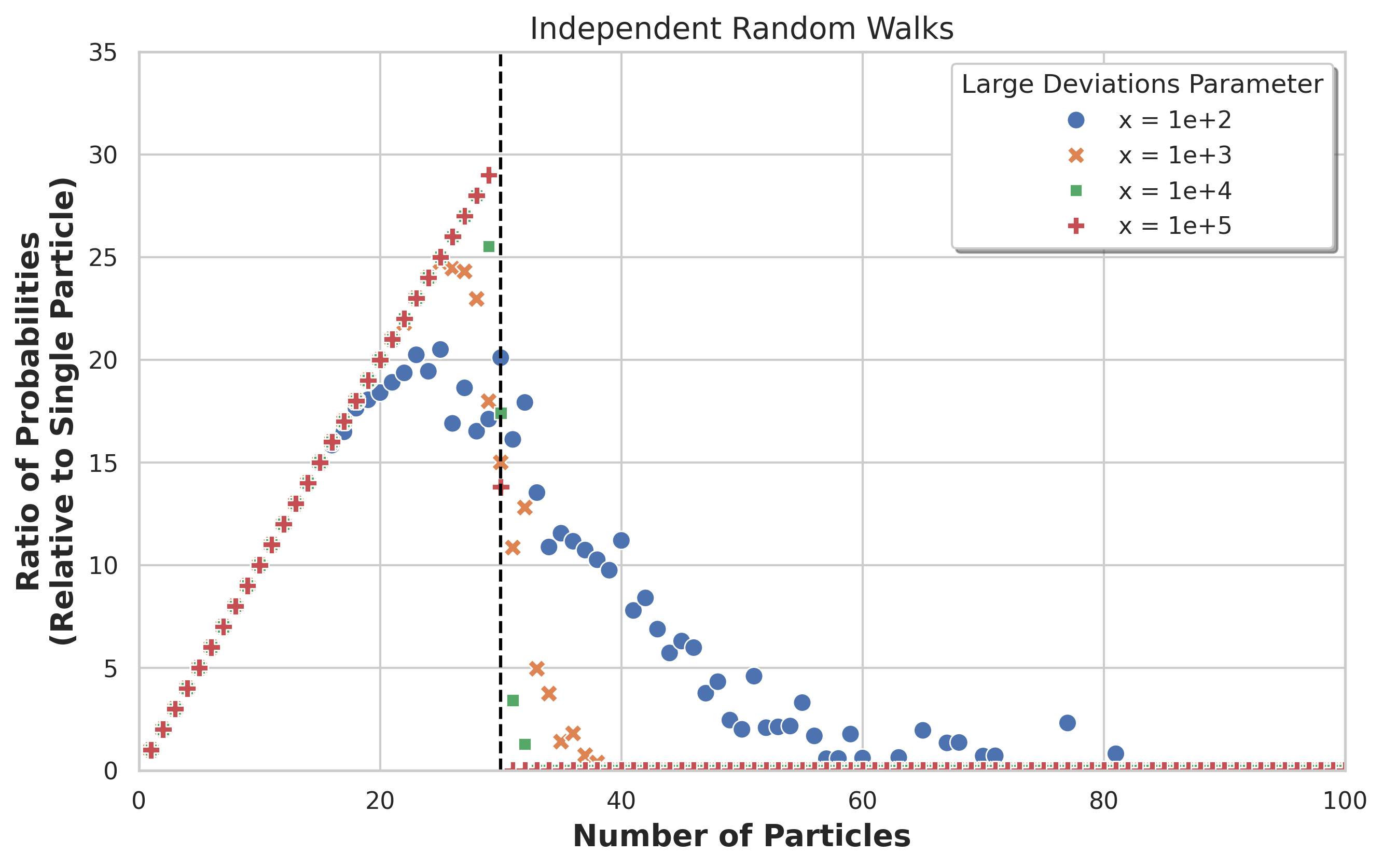}
    \caption{Parallel exploration for random walks with negative mean $p - (1-p)=-0.1$. Estimated ratio between $\PP(\tau^{(N)}(x)\leq B(x)/N)$ and  $\PP(\tau(x)\leq B(x))$ as a function of number $N$ with $B(x)=C\cdot x = 300\cdot x$ . The phase transition is observed at the expected threshold $N^* = \lceil C (1-2p)\rceil - 1=29$. }
    \label{fig:RW_Parallel_simulations}        
\end{figure}

Consider the random walk case described in Example \ref{example:BD_chains}. We use the parameters $p=0.45$ and $B(x)=300\cdot x$, which imply $\psi'(\lambda^*)=0.1$.  According to Theorem \ref{thm:main_theorem_num_particles} the ratio
\begin{equation}\label{eq:quotient_parallel_exp}
    \frac{\PP\left(\tau^{(N)}(x)\leq \frac{B(x)}{N}\right)}{\PP\left(\tau(x)\leq B(x)\right)}
\end{equation} is expected to approach the identity for $N< 300\cdot 0.1=30$ and zero for $N>30$, as $x\to\infty$. This is the behaviour observed in Figure \ref{fig:RW_Parallel_simulations}.

For the simulation we used the measure $\PP^{\lambda^*}$ as defined in \eqref{eq:exp_change_of_measure},  under which the process is distributed as a birth-and-death random walk with interchanged parameters, namely $p^* = \PP^{\lambda^*}(X_1 = 1) =1-\PP^{\lambda^*}(X_1 = -1) = 1-p$. Under this measure, passage over a given $x\in\NN$ is not a rare event anymore so simulations are feasible. For each large deviation parameter $x\in\{20,100,500,1000,2500\}$, and for each number of particles $N$ between $1$ and $100$ we simulate and average over $1000$ random walks to get an estimate of $\PP^{\lambda^*}(\tau \leq B(x)/N)$. By reverting the measure change we obtain estimates for the original ratio. %{\color{red} SHOULD WE SPECIFY THE FOLLOWING?  We then plot the quotient $N\frac{\PP^{\lambda^*}(\tau \leq B(x)/N)}{\PP^{\lambda^*}(\tau \leq B(x))}$ which approximates \eqref{eq:quotient_parallel_exp} up to an error of order at most $C_N\exp(-\lambda^* x)$ (see equation \eqref{eq:minimum_of_passage_times} below). }

\subsection{Parallel exploration driven by a L\'evy process}

Random walks have an analogue in continuous time known as Lévy processes. Under analogous assumptions on the exponential moments, we can extend Theorem \ref{thm:main_theorem_num_particles} to the case of Lévy processes on $\mathbb{R}$. This is based on the fact the asymptotics for first passage times over high barriers are similar in the cases of random walks and Lévy processes, see Section \ref{sec:proofs_parallel}. We now review some preliminaries on L\'evy processes, all of which may be found in classical references such as \cite{Levy_processes-Kyprianou} and \cite{Levy_Processes-Bertoin}.

%\seva{We need to think about the journal we are submitting to and decide whether we need to keep the detailed definitions.}

Let $(\Omega,\cF,\{\cF_t\}_t,\PP)$ be a filtered probability space, where the filtration $\{\cF_t\}_t$ is assumed augmented and right-continuous. A L\'evy process $Z = \{Z(t)\}_{t\in I}$, $I=[0,+\infty)$ is a stochastic process such that for all $0\leq s\leq t$, $Z(t)$ is $\cF_t$ measurable, $Z(0)=0$, $Z(t)-Z(s)$ is independent of $\cF_s$ and has the same distribution as $Z(t-s)$, and is continuous in probability, namely that $Z(t+s)\to Z(t)$ in probability if $s\to 0$. We will always work with a c\'adl\`ag version of $Z$. Let $\psi$ denote the L\'evy exponent of $Z$: \begin{equation*}
    \EE e^{\lambda Z(t)} = e^{t\psi(\lambda)}.
\end{equation*} The L\'evy-Khintchine formula introduces the characteristic triplet $(-\mu,\sigma,\Pi)$ of $Z$:
\begin{equation}\label{eq:levy_exponent}
    \psi(\lambda) = \log \EE e^{\lambda Z(1)} = -\mu \lambda + \frac{(\sigma\lambda)^2}{2} + \int_{-\infty}^{+\infty} (\exp(\lambda y)-1 - \lambda y \bfone(|y|\leq 1))\Pi(dy);
\end{equation} where $\mu\in\RR$ is called the drift coefficient, $\sigma\in (0,+\infty)$ is the diffusion coefficient (note that we exclude the finite-variation case) and $\Pi$ is a measure in $\RR\setminus\{0\}$ such that \begin{equation*}
    \int_{-\infty}^{+\infty} (1\wedge y^2) \Pi(dy) <+\infty,
\end{equation*} called the L\'evy measure (or jump measure). As in the discrete-time case we will work under Cram\'er's condition: $\EE(e^{\lambda Z(1)})$ is finite for $\lambda\in [0,\lambda_{\max})$ and there exists $\lambda^*\in (0,\lambda_{\max})$ such that $
    \psi(\lambda^*) = 0$.
This condition implies in particular that $\EE Z(1)< 0$ and that $Z$ drifts to $-\infty$ almost surely.

Let $\tau(x)$ and $\tau^{(N)}(x)$ be defined as \eqref{def:passage_time} and \eqref{def:min_of_passage_times} respectively. Up to taking a continuous time parameter, our result on the number of particles under finite time budget constraints for L\'evy processes is the same as for random walks:

\begin{theorem}\label{thm:main_theorem_num_particles_Levy}
    Let $Z$ denote a L\'evy process such that $Z(1)$ satisfies assumptions \eqref{right_cramer_cond}, \eqref{diff_mgf} and \eqref{positive_cramer_exp}. Let $B(\cdot)$ belong to the class $L(\lambda)$ for some $\lambda\in\Lambda_+$ (as defined in \eqref{def:asymp_linear_funcs}). Then for $N\geq 2$:
    \begin{equation}\label{eq:threshold_number_particles_Levy}
        \lim_{x\to+\infty}\frac{\PP\left(\tau^{(N)}(x)\leq \frac{B(x)}{N}\right)}{\PP\left(\tau(x)\leq B(x)\right)} =  \begin{cases}
            N,\ \text{ if } N \psi'(\lambda) < \psi'(\lambda^*);\\
            0,\ \text{ if } N \psi'(\lambda)  > \psi'(\lambda^*).
        \end{cases}
    \end{equation}
\end{theorem}

In view of the above result we obtain an optimal number of particles for L\'evy process, in the same way as in Corollary \ref{cor:optimal_num_part}.

The following is an example where all the parameters can be computed explicitly.

\begin{example}\label{example:BM_dirft}\textbf{Linear Brownian motion:}
     If the L\'evy measure is zero everywhere, the process is a Brownian motion with drift, namely a solution of the SDE \[dZ(t) = -\mu dt + \sigma dB(t)\ . \] We assume further that $\sigma=1$. Since its distribution at time $1$ is Gaussian with parameters $(-\mu,1)$, from Example \ref{example:gaussian_increments_rw}  we conclude that $\lambda^*=2\mu$ and $\psi'(\lambda^*) = \mu$. The same conclusion on the optimal number of particles as in Example \ref{example:gaussian_increments_rw} holds.\rqed
\end{example}

\subsubsection{Simulations for parallel L\'evy processes:}\label{subsec:simulations_LP_parallel}

\begin{figure}
    \centering
    \includegraphics[width=0.9\linewidth]{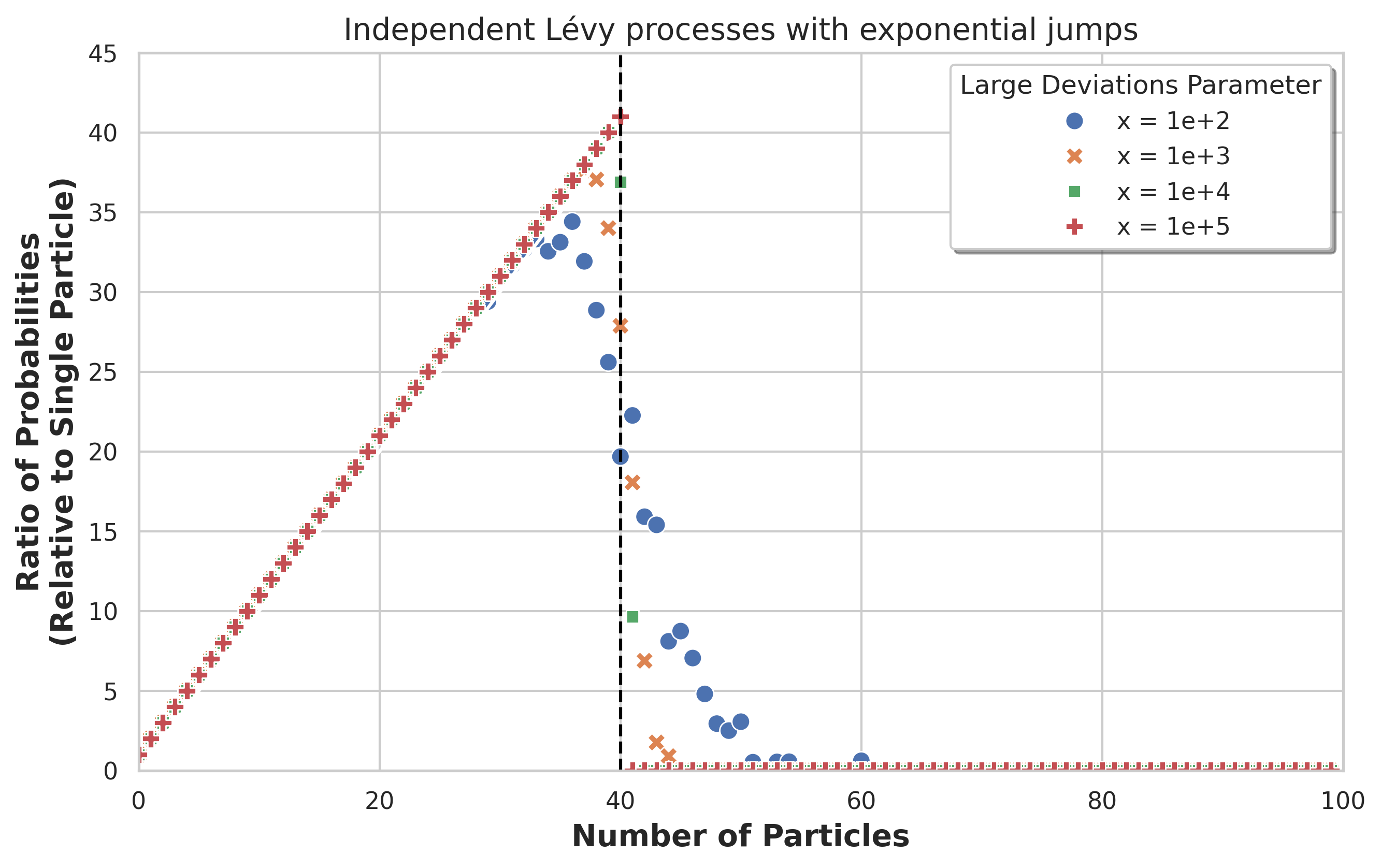}
    \caption{\textbf{Parallel exploration L\'evy processes}. Estimated ratio between $\PP(\tau^{(N)}(x)\leq B(x)/N)$ and $\PP(\tau(x)\leq B(x))$ as a function of the number $N$ of independent L\'evy processes with exponential jumps, with parameters described in Section \ref{subsec:simulations_LP_parallel}.}
    \label{fig:LP_parallel_simulations}
\end{figure}

For the simulations in Figure \ref{fig:LP_parallel_simulations} we considered a L\'evy process with positive and negative jumps. Positive (resp. negative) jumps occur at times distributed as a Poisson process of intensity $r=2$ (resp. $s=3$) and whose lengths are exponentially distributed with rate $\alpha=4$ (resp. $\beta = 1$). Following the notation in \eqref{eq:levy_exponent} we further take $\mu = \sigma = 1$. It can be checked that the L\'evy measure is given by: 
\begin{equation*}
    \Pi(dy) = \bfone(y< 0)3\exp(y)\; dy +  \bfone(y >0) 8\exp(-4y) \;dy.
\end{equation*}

 For the considered process, $\psi(\lambda)= -\lambda + \frac{1}{2}\lambda^{2}- \frac{3\lambda}{\lambda + 1}+ \frac{2\lambda}{4 - \lambda}, \;\; \lambda < 4$ and its positive root turns out to be $\lambda^* = 2$. Moreover, $ \psi'(2) = 8/3$ which gives, for $B$ taken as $B(x) = 15\cdot x$, an optimal number $N^* = 40$ of particles.

For simulation purposes, as in the discrete-time case, we use the exponential martingale associated to $\lambda^*$ and the measure $\PP^{\lambda^*}$, as defined in Section \ref{sec:proofs_parallel}. Specifically, the parameters of the process were chosen in such a way that under the measure $\PP^{\lambda^*}$, the process is again L\'evy with exponential jumps. Moreover the time intensities and the jump rates are permuted: under the measure $\PP^{\lambda^*}$ the parameters of the L\'evy measure are $r^*=4$, $\alpha^* = 2$, $s^* = 1$ and $\beta^* = 3$. 

\begin{remark}{\textbf{ A remark on reflected process}:} 
    In terms of applications it may be desirable to include reflection at a low barrier in our models to convey the idea that complexity of an exploration task cannot diverge to $-\infty$ in realistic scenarios. However, in the context of rare events we are in, there is no substantial difference in introducing reflection, since estimates vary at most by a multiplicative factor. To see this, consider the process $\Tilde{Z}(t) = Z(t)-[\inf_{0\leq s\leq t}Z(s)\wedge(-1)]$, namely the reflected version of $Z$ at $-1$. 
By considering the coupling $(Z(t),\Tilde{Z}(t))$ and observing that $\PP(Z(t)\leq\Tilde{Z}(t))=1$ for every $t\geq 0$, it follows that \begin{equation*}
        \PP(\tau(x)\leq t) \leq \PP(\Tilde{\tau}(x)\leq t), \text{ for every } t\geq 0. 
    \end{equation*}
    
    On the other hand, the event $\{\Tilde{\tau}(x)< t\}$ may be decomposed into trajectories that reach $x$ before the first reflection at $-1$ and those that do not. The first subset may be treated as a process without reflection, thus contributing a term of the same order as that of the process without reflection. Trajectories of the latter subset have a last reflection time before $\Tilde{\tau}(x)$, so there is an excursion of length $x+1$ during a shorter time. This provides an extra term of order at most $\PP(\tau(x+1)\leq t)<\PP(\tau(x)\leq t)$. Thus \[\PP(\Tilde{\tau}(x)\leq t) \leq 2\PP(\tau(x)\leq t)\] and we may then safely restrict ourselves to processes without reflection. The interested reader may consult \cite{reflected_Levy_Process_Maximum} for the maximum of reflected L\'evy processes.
    \rqed
    \end{remark}

%%%%%%%%%%%%%% Restart

\subsection{Exploration with restart for random walks and L\'evy processes}\label{subsec:Model_restart}

Our Strategy II consists of restarting a random walk or a L\'evy process as those studied in Theorems \ref{thm:main_theorem_num_particles} and \ref{thm:main_theorem_num_particles_Levy} upon leaving $(0,x)$ with a given probability measure $\nu_x$ supported in the interval $(0,x)$. Our first result in this direction, stated in Theorem \ref{thm:restart_main_result_general_measure}, gives asymptotic performance guarantees of restart strategies under mild assumptions on the restart measures $\nu_x$, as $x\to+\infty$.

In the next section we will specialize to a particularly interesting family of measures: quasi-stationary measures, for which better estimates are obtained (in particular matching bounds up to multiplicative constants). We remark that for the QSD case we rely exclusively on quasi-stationarity and related properties but not on Theorem \ref{thm:restart_main_result_general_measure}.

\begin{definition}
    Let $\nu_1,\nu_2$ be two positive finite measures on $\RR$. $\nu_1$ is said to be \textit{stochastically dominated} by $\nu_2$, denoted $\nu_1\preceq_{\text{stoch}} \nu_2$, if for every measurable non-decreasing function $u$: \begin{equation}\label{eq:sothc_domination_def}
        \int_{\RR} u(y)\ \nu_1(dy) \leq \int_{\RR} u(y)\ \nu_2(dy).
    \end{equation}
\end{definition} 

This notion is sometimes called \textit{first-order stochastic domination} and is equivalent to $\nu_1(x,+\infty)$ being less or equal than $\nu_2(x,+\infty)$ for every $x\in\RR$. 

The next theorem provides estimates on the time it takes for the restarted process to reach the desired region (here an interval $[x,+\infty)$ for large $x$) under mild assumptions on the restart measures $(\nu_x)_{x\geq 1}$, for which we assume without loss of genearality that $x\geq 1$.

\begin{theorem}\label{thm:restart_main_result_general_measure}
Let $\{Z(t)\}_{t>0}$ denote a random walk or a L\'evy process satisfying Cram\'ers condition, \textbf{and let $\lambda^*$ be given by \eqref{positive_cramer_exp}}. For each $x>0$ let $\{Z^{\nu_x}(t)\}_{t>0}$ be the restarted process with restart measures $(\nu_x)_{x\geq 1}$ as described at the beginning of this section, and let $\tau^{\nu_x}(x)$ be the first time a cycle of the process ends above $x$. Assume that there exists a finite positive measure $\nu$ with the following properties: \begin{enumerate}
    \item $\nu$ is not a delta measure on $0$,
    \item $\nu_x \stochdomd \nu$ for every positive $x$,
    \item $\nu$ has a finite second moment: $\int y^2\ \nu(dy)<+\infty$.
\end{enumerate}
Then, $B(x)$ growing faster than $\psi'(\lambda^*)^{-1}\;x$ but slower than $e^{\lambda^* x}$, the ratio 
\begin{equation}\label{eq:restart_main_result}
    \frac{\PP(\tau^{\nu_x}(x)\leq B(x))}{\PP(\tau(x)\leq B(x))} \frac{1}{B(x)\int_0^x\exp(\lambda^* y)\nu_x(dy)}
\end{equation} is bounded away from zero uniformly in $x>0$.
\end{theorem}

A typical trajectory of the restarted process is a concatenation of independent cycles, each of which consists of a trajectory $\{y_0 + Z(t): 0\leq t \leq \tau\}$ for some $y_0\in (0,x)$ sampled according to $\nu_x$ and $\tau \coloneqq \inf\{t>0: y_0 + Z(t) \notin (0,x)\}$.  An intuition behind the result is that adding $\nu_x$ as an initial distribution boosts the performance by a factor $\int_0^x\exp(\lambda^* y)\nu_x(dy)$ for each cycle. Since the cycles' duration is stochastically bounded we expect to observe an asymptotically-linear amount of cycles within time $B(x)$. However a matching upper bound remains elusive for this general context.

\begin{example}\label{example:1_restart_measure}
    Consider some fixed probability measure $\nu$ on $(0,+\infty)$ such that $\nu((0,1))>0$ and satisfying properties $1,2$ and $3$ listed in Theorem \ref{thm:restart_main_result_general_measure}. Now define the sequence $(\nu_x)_{x\geq 1}$ as the conditional laws \[\nu_x(\;\cdot\; ) \coloneq \frac{\nu(\;\cdot\;\cap\; (0,x))}{\nu((0,x))}. \] Then $\nu_x(y,+\infty) \leq \nu((0,1))^{-1}\; \nu(y,+\infty)$ and Theorem \ref{thm:restart_main_result_general_measure} is applicable for this family of measures.
\end{example}

\begin{remark}\label{remark:weak_convergence_and_domination}
    In a more general fashion than the previous example, one may wonder if the hypothesis of Theorem \ref{thm:restart_main_result_general_measure} are satisfied whenever the measures $\nu_x$ have a weak limit $\nu$. This is indeed the case under some assumptions on the moments of the $\nu_x$'s, but the dominating measure need not be related to the weak limit. Define $\nu_*((-\infty,y])\coloneq \inf_{x\geq 1} \nu_x((-\infty,y])$. Then $\nu_*$ defines a càdlàg function with $\lim_{y\to-\infty} \nu_*((-\infty,y]) = 0$. Since the $\nu_x$'s have a weak limit, $\lim_{y\to\infty} \nu_*((-\infty,y]) = 1$ also holds and $\nu_*$ defines a probability measure that dominates $\nu_x$ for every $x\geq 1$. Note as well that $\nu_*$ is not a delta measure at zero unless $\nu_x =\delta_0$ for all $x\geq 1$. Let us now give a precise condition under which $\nu_*$ satisfies property $3$ in Theorem \ref{thm:restart_main_result_general_measure}. Assume that the $\nu_x$'s have uniformly bounded moments of order $p$ for some $p>2$, i.e.: \begin{equation} \label{eq:moment_assumption}
       C\coloneqq \sup_{x\geq 1}\int_0^{+\infty} y^p \;\nu_x(dy) <+\infty.
    \end{equation} We bound the tails of the $\nu_x$ by Markov's inequality:\[\nu_x(y_0,+\infty) \leq y_0^{-p} \int_0^{+\infty} y^p\; \nu_x(dy) \leq  C y_0^{-p}. \]
    Then using the Layer Cake Representation the desired property holds: 
    \begin{align*}
        \int_0^{+\infty} y^2 \;\nu_*(dy) = 2 \int_0^{+\infty} y \;\nu_*(y,+\infty)\; dy &= 2 \int_0^{+\infty} y \; \sup_{x\geq 1}\nu_x(y,+\infty)\; dy \\ &\leq 2C\int_0^{+\infty} y^{1 - p} dy  < +\infty.
    \end{align*}
\end{remark}
As a final example for this section we take a simple extension of Example \ref{example:1_restart_measure} that shows that Theorem \ref{thm:restart_main_result_general_measure} covers cases beyond weak convergence:

\begin{example}
Take $\nu^1$ and $\nu^2$ two (distinct) probability measures on $(0,+\infty)$ satisfying the conditions of Theorem \ref{thm:restart_main_result_general_measure} and $\nu^i((0,1))>0,\ i=1,2$. We may build an alternating family of measures for $x\geq 1$ as follows: 

\[\nu_x(\;\cdot\; ) \coloneq 
\begin{cases}
  \frac{\nu^1(\;\cdot\;\cap\; (0,x))}{\nu^1((0,x))} \text{ if } \lfloor x \rfloor\text{ is odd, }\\ 
  \frac{\nu^2(\;\cdot\;\cap \;(0,x))}{\nu^2((0,x))} \text{ if } \lfloor x \rfloor\text{ is even. }
\end{cases}  
\]  In this case Theorem \ref{thm:restart_main_result_general_measure} applies, but there is not a weak limit of the whole sequence of measures.
\end{example}

%%%%%%%%%%%%%%%%%%%%%%%%%%%%%%%% Restart QSD

\subsection{Restart with quasi-stationary measures.}\label{subsec:Model_restart_qsd}

We now specialize to an important family of restarting measures for which we can obtain sharper asymptotic results than those provided by Theorem~\ref{thm:restart_main_result_general_measure}.  

Throughout this section we assume that $\nu_x$ is the quasi-stationary distribution (QSD) of the process restricted to the interval $(0,x)$.  
We begin with a few brief preliminaries on quasi-stationary distributions, and refer the reader to the monograph \cite{QSD_book_Collet_etal} for a comprehensive treatment of the subject.

\begin{definition}[Quasi-stationary measure]
    Let $Z=\{Z(t)\}_{t\geq 0}$ be a Markov process with state space $\cS$ and let $A$ be a subset of $\cS$. The absorbed process, denoted by $Z_{A}$, is defined as $Z_{A}(t) = Z(t\wedge \tau(A))$, where $\tau(A) \coloneqq \inf\{t>0: Z(t)\in A \}$. A quasi-stationary measure (QSD) $\nu$ is a probability measure on $\cS\setminus A$ which is invariant when conditioned on non-absorption, which means: \[\PP_{\nu}(Z(t)\in B\ |\ t<\tau(A) ) = \frac{\PP_{\nu}(Z_A(t)\in B)}{\PP_{\nu}(t<\tau(A) )} = \nu(B),\] for every $t>0$ and every measurable set $B\subset \cS\setminus A$. 
\end{definition}

The study of QSD for a given Markov process is a subtle matter, in general they may not exist and when they do they may not be unique. 

 For the particular case of a one-dimensional Lévy process with absorbing set $A = \RR\setminus (0,a)$ for some $a>0$, by \cite[Theorem 2.1 and Remark 2.2]{Kolb2014} there is a unique QSD, provided the distribution of the process without absorption at any time $t>0$ is absolutely continuous with respect to Lebesgue in $\RR$, which is satisfied \textbf{since} we are assuming a strictly positive diffusion \textbf{coefficient. 
%at any time $t>0$ a L\'evy process is a sum of a deterministic (drift) process, a pure-jump process and an independent Brownian motion. Hence, the distribution of the process at any given time is a convolution involving a Gaussian factor, which implies the absolute continuity of the distribution of the process. 
When} there is a unique QSD, it is obtained as \textit{the Yaglom limit} \[\lim_{t\to+\infty} \PP_y(Z(t)\in \cdot\ |\ t < \tau(A));\] which is independent of the starting position $y\in \cS\setminus A$.

Existence of a Yaglom limit for a L\'evy process when the absorbing set is $(-\infty,0]$ was proved in \cite{Kyprianou_Palmowski_QSD_LEVY}\footnote{The processes we consider fall in their class A category with parameter $\alpha=2$.}, extending the case without jumps of \cite{martinez_san_martin_QSD-BM_drift}.

To the best of our knowledge, there is no analogous existence-and-uniqueness result specific to discrete-time random walks conditioned on non-absorption on bounded intervals. If the law of the increments of a random walk is absolutely continuous with respect to the Lebesgue measure, then one may consider the associated compound Poisson process (CPP) with rate $1$, thereby obtaining a L\'evy process with an absolutely continuous law for every time larger than its first jump. By checking that for large $t$ the probability of not having jumped is negligible for the CPP \textit{even when conditioned on non-absorption by time $t$}, one obtains the same Yaglom limit as the conditioned CPP and in particular existence and uniqueness given by \cite{Kolb2014}. This can be carried out using the asymptotics on passage times from \cite{Seva_passage_times_LP_RW}. These asymptotics for conditioned processes are used in the proof of Lemma \ref{lemma:qsd_converge_weakly} of Section \ref{sec:proof_restart}.
Yaglom limits for random walks with absorption at $(-\infty,0]$ were studied in \cite{Iglehart_yaglom_limit_RW}, and under Cram\'er condition in \cite{doney_cond_RW_cramer_condition}.

We may now state our main theorem for restarted processes:

\begin{theorem}\label{thm:main_thm_QSD_Levy}
Let $Z$ be a random walk on discrete time or a  L\'evy process satisfying \eqref{right_cramer_cond}, \eqref{diff_mgf} and \eqref{positive_cramer_exp} as in Theorem \ref{thm:main_theorem_num_particles}. For each positive $x\in (0,+\infty)$ let $\nu_x$ denote the quasi-stationary probability measure on $(0,x)$ for the process $Z$ absorbed upon exiting $(0,x)$. Let $Z^{\nu_x}$ be the restarted version of $Z$ on $(0,x)$ with \textbf{restart} measure $\nu_x$ and  $\tau^{\nu_x}(x)$ the associated passage time over $x$. Let the time budget $B(x)$ grow faster than $x/\psi'(\lambda^*)$ as $x\to+\infty$ but slower than $e^{\lambda^* x}$. Then there exist constants $c,C>0$ such that:
    \begin{equation}\label{eq:gain_restart_mech} 
        c \leq \frac{\PP(\tau^{\nu_x}(x)<B(x))}{\PP(\tau(x)<B(x))}\cdotp \frac{1}{B(x) \int_0^x e^{\lambda^* y}\nu_x(dy)} \leq C
    \end{equation}
    for all $x>0$.
\end{theorem}

\begin{corollary}\label{Coro:linear_BM_QSD}
    In the particular case when $Z$ is a linear Brownian motion, with a negative drift equal to $-\mu$, we have \begin{equation}
       \int_0^x e^{\lambda^* y}\nu_x(dy) =  e^{\mu x},
    \end{equation} and the ratio\begin{equation}\label{eq:corolary_BM_case}  
        \frac{\PP(\tau^{\nu_x}(x)<B(x))}{\PP(\tau(x)<B(x))} \cdotp \frac{1}{B(x) e^{\mu x}}
    \end{equation} is bounded away from zero and infinity for $B(x)$ growing faster than $x/\mu$. 
\end{corollary}
    
The theorem states that the application of a restarting mechanism with quasi-stationary distributions improves performance by a factor comparable with the time budget times an exponential moment of the respective measure. In the Brownian motion case, this improvement factor is exponential in the complexity of the task. Indeed, while $\PP(\tau(x)\leq B(x))$ is of order $e^{-2\mu x}$, with the restarting mechanism we get order $B(x)e^{-\mu x}$. For a general L\'evy process, exponential moments of its quasi-stationary measures are not easy to estimate. We conjecture that they grow as $\exp\left[(\lambda^*-\lambda_0) x\right]$, as $x\to+\infty$, where $\lambda_0$ is the positive solution of $\psi'(\lambda) = 0$. Note that this is exactly the case for the linear Brownian motion. For spectrally negative L\'evy processes, namely when $\Pi(0,+\infty) = 0$, the conjecture might be proved with the use of a Girsanov transformation of the process into a zero-drift one (thus providing the conjectured exponential factor) and applying the results of \cite{Lambert_qsd_finite_interval}. We do not pursue this further here, it is left as a direction of further study.

%%%%%%%%%%%%%

\subsubsection{Simulations for restarted L\'evy processes}\label{sec:simulations_restart}

We now present numerical experiments for restarted processes. As a case study,
we consider the Lévy process with exponential jumps introduced in
Subsection~\ref{subsec:simulations_LP_parallel}, equipped with a restart mechanism
activated upon exiting the interval $(0,50)\subset\mathbb{R}$.

For the process without restart, recall that the Cram\'er exponent is
$\lambda^* = 2$. Consequently, when the process is started from the origin,
the probability of ever reaching the upper barrier at $50$ is of order
$\exp(-100)$.
As restart measure, we consider an exponential distribution with mean $10$,
truncated to $(0,50)$, namely
\[
\nu_{50}(dx)
=
\frac{0.1\,e^{-0.1x}}{1-e^{-5}}\,\mathbf{1}_{(0,50)}(x)\,dx.
\]
Its exponential moment of order $\lambda^*$ satisfies
\[
\int_0^{50} e^{2x}\,\nu_{50}(dx)
\;\approx\;
9.6\times 10^{39}.
\]
Multiplying by $\mathbb{P}_0(\tau(50)<+\infty)$ yields
\[
\mathbb{P}_0(\tau(50)<+\infty)
\int_0^{50} e^{2x}\,\nu_{50}(dx)
\;\approx\;
3.6\times 10^{-4}.
\]

\begin{figure}
    \centering
    \includegraphics[width=.9\linewidth]{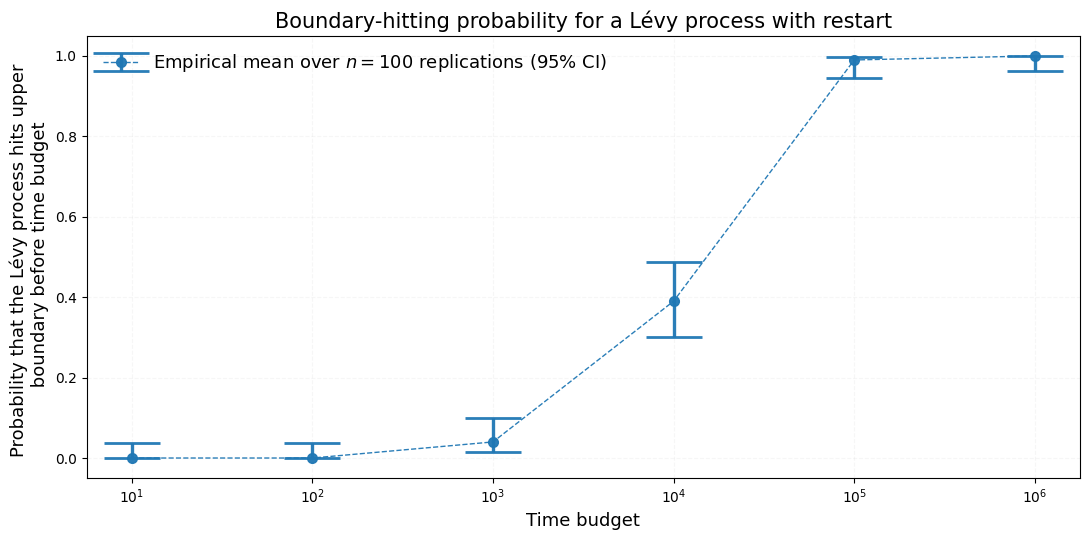}
    \caption{
        Probability of exceeding the level $x=50$ for the restarted Lévy process as a
        function of the time horizon.
        The underlying process combines a Brownian motion with drift $\mu=-1$ and
        volatility $\sigma=1$ with exponential jumps: positive jumps arrive at rate
        $\lambda_+=2$ with jump sizes of rate $4$, while negative jumps arrive at rate
        $\lambda_-=3$ with jump sizes of rate $1$.
        Whenever the process exits the interval $(0,50)$, it is restarted according to
        a truncated exponential distribution with mean $10$.
       }
    \label{fig:Restarted_Levy}
\end{figure}

Although this estimate is only heuristic, in view of
Theorem~\ref{thm:restart_main_result_general_measure} it suggests that the
probability of exceeding the level $x=50$ within a finite time horizon
$B(50)$ under this restart strategy is no longer negligible.
This behavior is confirmed numerically in
Figure~\ref{fig:Restarted_Levy}, which reports empirical mean estimates based on
$100$ independent replications for each time budget, together with
$95\%$ confidence intervals.  The figure illustrates how the restart mechanism makes the exceedance of the
        upper barrier observable on moderate time scales, with confidence intervals (CI)
        indicating the associated variability.

We emphasize that the simulations reported in this section do not rely on importance
sampling or on any other variance–reduction technique. The numerical results therefore illustrate that
introducing a restart mechanism alone can already lead to a significant
improvement in the simulation of rare events. In particular, restarting makes it
possible to turn probabilities that are extremely small for the original process
into quantities that can be reliably estimated within reasonable computational
time.

%%%%%%% SECTION PROOFS THMs 1 and 2

\section{Proof of Theorems \ref{thm:main_theorem_num_particles} and \ref{thm:main_theorem_num_particles_Levy}}\label{sec:proofs_parallel}

%{\color{purple} 

We begin the section by introducing some preliminaries on exponential martingales, first for random walks and then for L\'evy process, which will be used in the proofs. They are standard techniques for the study of large deviations events for light tails, as is our case. 

\subsection{Preliminaries}

By assumption \eqref{right_cramer_cond}, the family of exponential martingales indexed by \( \lambda \in \Lambda_+ \) may be defined as:
\[
    M^{\lambda}_t \coloneq \exp\left( \lambda Z(t) - t\psi(\lambda) \right), \quad t \in \mathbb{N}.
\]
Indeed, \( M^{\lambda} \) is a martingale for every \( \lambda \in \Lambda_+ \), which acts as a density for a new probability measure on the same probability space. Specifically, we define the measure \( \mathbb{P}^{\lambda} \) as
\begin{equation}\label{eq:exp_change_of_measure}
    \mathbb{P}^{\lambda}(A) \coloneq \mathbb{E}\left[M_{\tau}^{\lambda} \mathbf{1}(A)\right],
\end{equation}
where \( A \) depends on \( \{Z(t)\}_{t \leq \tau} \), and \( \tau \) is any almost-surely finite stopping time. The measure \( \mathbb{P}^{\lambda^*} \), associated with \( \lambda^* \) in equation \eqref{positive_cramer_exp}, will be of particular importance.

We also observe that the density of \( \mathbb{P}|_{\mathcal{F}_{\tau}} \) with respect to \( \mathbb{P}^{\lambda}|_{\mathcal{F}_{\tau}} \) is given by \( (M^{\lambda}_{\tau})^{-1} \). 
%Expectation with respect to \( \mathbb{P}^{\lambda} \) will always be denoted by \( \mathbb{E}_{\lambda} \).

Exponential martingales are particularly useful for computing probabilities of rare events, as these events may become more frequent under \( \mathbb{P}^{\lambda} \) if \( \lambda \) is appropriately chosen. Under the measure \( \mathbb{P}^{\lambda} \), the increments \( X_j \), \( j \geq 1 \), are i.i.d. (making \( Z \) a random walk), with mean and variance given by:
\[
\mu(\lambda) \coloneqq \mathbb{E}_{\lambda}[X_1]  = \psi'(\lambda) \quad \text { and } \quad
\sigma(\lambda) \coloneqq \mathbb{E}_{\lambda}\left[ |X_1 - \mu(\lambda)|^2 \right] = \psi''(\lambda).
\]
The threshold in Theorem \ref{thm:main_theorem_num_particles} can thus be interpreted in terms of the measure \( \mathbb{P}^{\lambda^*} \): the random time $\tau(x)$ is $\PP^{\lambda^*}$-almost surely finite (since $\psi'(\lambda^*)>0$) with mean \( x \psi'(\lambda^*)^{-1} \), asymptotically as \( x \to +\infty \).
%}

%{\color{blue} 
For L\'evy process we have analogous definitions: under Cram\'er's condition the family $\{M^{\lambda}_t\}_{0\leq t\leq \lambda_{\max}}$ of exponential martingales is defined as \begin{equation*}
    M^{\lambda}_t \coloneqq \exp(\lambda Z(t) - t\psi(\lambda)),
\end{equation*} and the corresponding family of measures is given by \begin{equation*}
    \frac{d\PP^{\lambda}}{d\PP}\Bigg|_{\cF_t} = M^{\lambda}_t.
\end{equation*}

Under any of the measures $\PP^{\lambda}$, $Z$ is again a L\'evy process with L\'evy exponent \begin{equation}\label{eq:change_measure_LP_exponent}
    \psi^{(\lambda)}(r) = \psi(\lambda + r)- \psi(\lambda),
\end{equation} with mean and variance again given by:
\[
\mu(\lambda) \coloneqq \mathbb{E}_{\lambda}[X_1]  = \psi'(\lambda)
\quad \text{ and } \quad 
\sigma(\lambda) \coloneqq \mathbb{E}_{\lambda}\left[ |X_1 - \mu(\lambda)|^2 \right] = \psi''(\lambda).
\]

From \eqref{eq:change_measure_LP_exponent} it is straightforward to check that the corresponding characteristic triplet is given by: 
\begin{equation}\label{eq:change_measure_LP}
    \mu_{\lambda} = -\mu+\lambda \sigma^2+\int_{-1}^1 y(e^{\lambda y }-1)\Pi(dy),\ \ \sigma_{\lambda}^2 = \sigma^2,\ \ \Pi^{\lambda}(dy) = e^{\lambda y}\Pi(dy). 
\end{equation}

The measure $\PP^{\lambda^*}$ provides an intuition to Theorem \ref{thm:main_theorem_num_particles_Levy} as in the random walk case: the law of the process conditioned to reach the high barrier is distributed as $\PP^{\lambda^*}$. We do not prove this exact result here but lies at the core of the references cited for our proof \cite{Cramer_estimate_Levy_process,Palmowski_etal_Levy_processes}.

\begin{remark}
    $\mu_{\lambda}$ in \eqref{eq:change_measure_LP} should not be confused with $\mu(\lambda)$. If the jump measure $\Pi$ is non-zero they might not coincide.\rqed
\end{remark}

%}

\subsection{Proofs}\label{subsect:Proofs}

\begin{remark}\label{rmk:asymp_notation}
 We introduce some pieces of notation: $f(x)\sim g(x)$ will be used to denote asymptotic equivalence $\lim_{x\to +\infty} f(x)g(x)^{-1} = 1$ and $f(x)\asymp g(x)$ for $ C^{-1}<f(x)g(x)^{-1}<C$ for some constant $C>0$ for all $x>0$.
\end{remark}

Theorems \ref{thm:main_theorem_num_particles} and \ref{thm:main_theorem_num_particles_Levy} will follow from the following asymptotics of the passage times of random walks and L\'evy processes with the assumptions on their exponential moments stated in Subsection \ref{subsec:Model_parallel}. 

The following proposition unifies \cite[Corollary 2.2]{passage_times_rw_Hoglund} for the case of random walks and \cite[Theorem 1]{Palmowski_etal_Levy_processes} for L\'evy processes:

\begin{proposition}
    \label{thm:estimates_pasasge_times}
    Let $Z = \{Z(t)\}_{t\in I}$ be either a random walk ($I=\NN$) or a Lévy process ($I=\RR^+$) taking values in $\RR$ and assume that conditions \eqref{right_cramer_cond}, \eqref{diff_mgf} and \eqref{positive_cramer_exp} are satisfied (for $Z(1)$ in the case of a L\'evy process). Fix some $\lambda\in (0,+\infty)\cap \Lambda$ such that $\mu(\lambda)=\psi'(\lambda)>0$. If $x$ and $t=t(x)$ go to infinity in such a way that $t\in L(\lambda)$
    %$x = \mu(\lambda) t + \mathbf{o}(\sqrt{t})$ 
    then there exist positive constants $C$ (independent of $\lambda$) and $D(\lambda)$ such that: \begin{equation}\label{eq:threshold_passage_times}
        \PP(\tau(x)\leq t) \sim \begin{cases}
            C \exp(-\lambda^* x), \ \text{if }\ \mu(\lambda) < \mu(\lambda^*)\\
            D(\lambda) t^{-1/2}\exp(-\zeta[\mu(\lambda)]t), \ \text{if }\ \mu(\lambda)>\mu(\lambda^*);
        \end{cases}
    \end{equation} where $\zeta$ is the convex conjugate of $\psi$: \begin{equation}\label{eq:convex_transform}
        \zeta[s] := \sup_{\lambda\in\Lambda}\{\lambda s - \psi(\lambda)\}.
    \end{equation}
\end{proposition}
%Recall that $\mu(\lambda)\coloneqq \psi'(\lambda) = \EE^{\lambda}Z(1)$. S
Some remarks are now in order: 
\begin{remark}
    The constants $C$ and $D(\lambda)$ are given in \cite{passage_times_rw_Hoglund} and \cite{Palmowski_etal_Levy_processes} in terms of the increments of the process in each case. Their explicit values are not relevant for our analysis.\rqed
\end{remark}
\begin{remark}\label{rk:rate_function_mu}
    If $s= \mu(\lambda)$, then $\lambda$ realizes the supremum in the definition of $\zeta[s]$, that is \begin{equation}
        \zeta[\mu(\lambda)] = \lambda \mu(\lambda) - \psi(\lambda),\ \lambda\in\Lambda.
    \end{equation}\rqed
\end{remark}

\begin{remark}\label{remark:oevrshoot_convergence}
    The case $\mu(\lambda) <\mu(\lambda^*) $ in the proposition shows that for $t$ growing sufficiently fast with $x$, $\PP(\tau(x)\leq t)$ and $\PP(\tau(x)< +\infty)$ have the same exponential order:
\begin{multline*}
    \PP(\tau(x)\leq t) \; \leq & \ \; \PP(\tau(x)<+\infty) \\ = &\ \; \EE_{\lambda^*}[\textbf{1}(\tau(x)<+\infty)\exp(-\lambda^*Z(\tau(x)))] \\= & \ \; e^{-\lambda^* x}\EE_{\lambda^*}[\textbf{1}(\tau(x)<+\infty)\exp(-\lambda^*[Z(\tau(x))-x])]\\ \leq &\ \; e^{-\lambda^* x};
\end{multline*}  
where we used that $\PP^{\lambda^*}(\tau(x)<+\infty) = 1$ and that the \textit{overshoot at $x$}, $U(x)\coloneq Z(\tau(x))-x$ is positive. In fact, if one assumes that the jump measure of $Z$ is not supported on a lattice then $\PP^{\lambda^*}(\tau(x)<+\infty)  \exp(\lambda^* x)$ converges to a positive constant (see \cite[Theorem 7.6]{fluctuations_LP_kyprianou}, \cite{Cramer_estimate_Levy_process},  \cite{Palmowski_etal_Levy_processes}). We will recall this fact repeatedly.\rqed
\end{remark}

% \begin{remark}
% Proposition \ref{thm:estimates_pasasge_times} provides the asymptotics to compare the probability of a particle reaching a given level $x$ within time $B(x)$, with that of the maximum of $N$ independent particles reaching the same level within time $B(x)/N$. As a consequence, an asymptotically optimal number of particles will be determined. As will be shown below, wether or not $N$ particles perform better than a single one for any given $N$, strongly depends on the side of the threshold on which $B(x)/N$ is located. \rqed
% \end{remark}

\begin{proof}{Proof of Theorems \ref{thm:main_theorem_num_particles} and \ref{thm:main_theorem_num_particles_Levy}:}

We begin by noticing that for any time $t$ (discrete or continuous) \begin{equation}\label{eq:minimum_of_passage_times}
    \PP\left(\tau^{(N)}(x)\leq t \right) = N \PP\left(\tau(x)\leq t \right) + \textbf{O}\left(\binom{N}{\lceil N/2\rceil} \PP\left(\tau(x)\leq t \right)^2\right),
\end{equation}
so for $N$ fixed, $\PP\left(\tau^{(N)}(x)\leq t \right)\sim N \PP\left(\tau(x)\leq t \right)$. Recall that $\lceil \cdot \rceil$ denotes the ceiling function.

  \textit{Case I: } $N\mu(\lambda) < \mu(\lambda^*)$.
    In this case $B(x)$ and $B(x)/N$ are in the same regime of Proposition \ref{thm:estimates_pasasge_times}, thus, by \eqref{eq:minimum_of_passage_times}: 
    \begin{equation}
       \lim_{x\to+\infty}\frac{\PP\left(\tau^{(N)}(x)\leq \frac{B(x)}{N}\right)}{\PP\left(\tau(x)\leq B(x)\right)}=  \lim_{x\to+\infty}\frac{N\PP\left(\tau(x)\leq \frac{B(x)}{N}\right)}{\PP\left(\tau(x)\leq B(x)\right)} = N.
    \end{equation}

For the remaining cases some considerations will be needed: for each $\lambda$ such that $\mu(\lambda)>0$ let 
\begin{equation}
S(\lambda) := \{s\geq 1: \exists \; \lambda_s>\lambda, \text{ s.t. } \mu(\lambda_s)=s\mu(\lambda)\}.    
\end{equation} The set $S(\lambda)$ contains a (non-trivial) right neighbourhood of $1$ as long as $\lambda < \lambda_{\max}$. We also note that for each $\lambda$ as before, the map $s\mapsto \lambda_s$ is continuously differentiable and \[\frac{d}{ds}\lambda_s = \frac{\mu(\lambda)}{\mu'(\lambda_s)} = \frac{\mu(\lambda)}{\sigma^2(\lambda_s)}.\]

 \textit{Case II: } $\mu(\lambda) <\mu(\lambda^*)<N\mu(\lambda)$. As a subcase, assume first that $N\in S(\lambda)$, so that there exists $\lambda_N\in (\lambda,\lambda_{\max})$ such that:
 \begin{equation}\label{eq:temporary_assumption_main_thm}
     \mu(\lambda_N) = N\mu(\lambda). 
 \end{equation}
    In this case, we have that $x = \mu(\lambda_N)\frac{B(x)}{N} + \mathbf{o}(\sqrt{x})$ so the exponential rate of $\PP\left(\tau(x)\leq \frac{B(x)}{N}\right)$ is given by Proposition \ref{thm:estimates_pasasge_times}:
    \begin{equation}
        \frac{B(x)}{N}\zeta[\mu(\lambda_N)] = \frac{B(x)}{N}(\lambda_N \mu(\lambda_N) - \psi(\lambda_N)).
    \end{equation}
    On the other hand $\PP\left(\tau(x)\leq B(x)\right)\sim e^{-\lambda^*x}$. Then, for some constant $\Tilde{D}(\lambda)$: 
    \begin{equation}\label{eq:case2_main_thm}
        \begin{aligned}
            \lim_{x\to+\infty}\frac{\PP\left(\tau(x)\leq \frac{B(x)}{N}\right)}{\PP\left(\tau(x)\leq B(x)\right)} & \leq \lim_{x\to+\infty}\Tilde{D}(\lambda) \exp\left( - x(\lambda_N-\lambda^*) + \frac{B(x)}{N}\psi(\lambda_N) \right) \\
            = \lim_{x\to+\infty}&\Tilde{D}(\lambda) \exp\left( - x(\lambda_N-\lambda^*) \left( 1 -\frac{\psi(\lambda_N)-\psi(\lambda^*)}{\mu(\lambda_N)(\lambda_N-\lambda^*)}\right)+\lambda_N\textbf{o}(\sqrt{x})\right) .
        \end{aligned}
    \end{equation}
    
    Using the Mean Value Theorem and the fact that $\mu$ is strictly increasing in the interval $(\lambda^*,\lambda_N)$ we conclude that the limit in \eqref{eq:case2_main_thm} is zero. 

For the second subcase, let $N\notin S(\lambda)$, so that $N\mu(\lambda)$ is at a positive distance above the range of $\mu(\cdot)$. A bound may be obtained from any $s\in S(\lambda)\setminus{\{1\}}$ by monotonicity of $t\mapsto \PP(\tau(x)\leq t)$.
\begin{equation}    
\begin{aligned}
         \lim_{x\to+\infty}\frac{\PP\left(\tau^{(N)}(x)\leq \frac{B(x)}{N}\right)}{\PP\left(\tau(x)\leq B(x)\right)} \leq      \lim_{x\to+\infty}\frac{\PP\left(\tau^{(N)}(x)\leq \frac{x}{\mu(\lambda_s)} + \mathbf{o}(\sqrt{x})\right)}{\PP\left(\tau(x)\leq B(x)\right)} = 0.
\end{aligned}
\end{equation}

\textit{Case III: } $\mu(\lambda^*)<\mu(\lambda)$. As in the subcases considered above, assume first that $N\in S(\lambda)$. Now both $B(x)$ and $B(x)/N$ lie in the right-hand side of the threshold in \eqref{eq:threshold_passage_times}. It suffices to show that the function 
\begin{equation*}
    s\in S(\lambda) \mapsto  \frac{1}{s}\zeta[\mu(\lambda_s)]
\end{equation*} is strictly increasing. Indeed, in that case \begin{equation}
  \lim_{x\to+\infty}\frac{\PP\left(\tau^{(N)}(x)\leq \frac{B(x)}{N}\right)}{\PP\left(\tau(x)\leq B(x)\right)} = \lim_{x\to+\infty}\exp\left\{ B(x)\left(\zeta[\mu(\lambda)] - \frac{1}{N}\zeta[\mu(\lambda_N)]\right)\right\}=0,
\end{equation} and the result is proved.

It only remains to prove that $\frac{1}{s}\zeta[\mu(\lambda_s)]$ is increasing; and using Remark \ref{rk:rate_function_mu} this is seen to hold: \begin{equation}
    \begin{aligned}
    \frac{d}{ds}\frac{1}{s}\zeta[\mu(\lambda_s)] & = \frac{d}{ds}\lambda_s \mu(\lambda) + \frac{\psi(\lambda_s)}{s^2} - \frac{1}{s} \mu(\lambda_s) \frac{d}{ds}\lambda_s \\& =  \frac{\psi(\lambda_s)}{s^2}, 
\end{aligned}
\end{equation} which is positive since $\lambda_s >\lambda> \lambda^*$.

The proof is finished by arguing in the same way as in Case II for $N\notin S(\lambda)$.$\bqed$
\end{proof}

%%%% Section 

\section{Proof of Theorem \ref{thm:restart_main_result_general_measure}}\label{sec:proof_restart_general_measure}

For the proof of Theorem \ref{thm:restart_main_result_general_measure} we will need the following auxiliary lemmas:

  %%%%%%%%%

\begin{lemma}\label{lemma:ruin_proba_order}
    Let
    \begin{equation}\label{eq:definition_of_q}
     q(x)\coloneqq \PP_{\nu_x}({\tau}(x)<{\tau}(0)).
    \end{equation} Then under the assumptions of Theorem \ref{thm:restart_main_result_general_measure}:   
    \begin{equation}\label{eq:q_to_0}
        \lim_{x\to+\infty} q(x) = 0.
    \end{equation} Moreover, there exists a constant $c>0$ such that: 
    \begin{equation}\label{eq:lema_3_statement}
      c \int_0^x e^{\lambda^* y}\nu_x(dy) \leq e^{\lambda^* x} q(x)  \leq \int_0^x e^{\lambda^* y}\nu_x(dy).
    \end{equation} 
\end{lemma}

    %%%%

  %%%%%%%%%%%%%%

\begin{lemma}\label{lemma:asymptotics_tau_0_tau_x}
    Let $\eta$ and $\zeta$ be random variables defined by: \[\PP(\eta \leq t) = \PP_{\nu_x} (\tau(0)\leq t\ |\ \tau(0)<\tau(x))\] and \[\PP(\zeta \leq t) = \PP_{\nu_x} (\tau(x)\leq t\ |\ \tau(x)<\tau(0)).\]
    Under the assumptions of Theorem \ref{thm:restart_main_result_general_measure}: \[\limsup_{x\to+\infty} \EE[\eta] <+\infty \] and\footnote{The notation $\asymp$ was introduced in Remark \ref{rmk:asymp_notation}.} \[\EE [\zeta]\ \asymp \frac{x}{\psi'(\lambda^*)}.\]
\end{lemma}
    
    %%%%%%%%%%%%%%

    %%%%%%%%%%%%%%%%%%%%%%%%
\begin{lemma}\label{lemma:E_delta_F_delta}
 Let $\eta,\zeta$ be as in Lemma \ref{lemma:asymptotics_tau_0_tau_x} and consider an i.i.d. sequence $\eta_1,\eta_2,\dots$ with the same distribution as $\eta$. For a fixed $\delta>0$ define the events \[E_{\delta}^m \coloneqq \left\{\sum_{j=1}^m \eta_j \leq m(1+\delta)\;\EE\eta_1\right\}\; \text{ and }\;F_{\delta} \coloneqq \bigg\{\zeta \leq (1+ \delta)\; \EE\zeta \bigg\},\] and for any $T>0$ let \[n_{\delta}(T) \coloneqq \bigg\lfloor\frac{T-(1+ \delta)\;\EE\zeta}{(1+ \delta)\;\EE\eta_1}\bigg\rfloor.\]
In the setting of Theorem \ref{thm:restart_main_result_general_measure} the following holds:  \begin{equation}\label{eq:E_delta_F_delta}
    \lim_{x\to+\infty}\PP\left( E_{\delta}^{n_{\delta}(B(x))}\ \cap\ F_{\delta}\right) = 1.
\end{equation}
\end{lemma}
   %%%%%%%%%%%%%%%%%%%%%%%%

Assuming these we may now give the main result for restarted processes.

\begin{proof}{Proof of Theorem \ref{thm:restart_main_result_general_measure}:}
    The first time the restarted process goes above $x$ has a regenerative structure with cycles determined by the restart times. Formally:
    \begin{equation}\label{eq:sum_decomposition_passage_time}
        \tau^{\nu_x}(x) = \sum_{j=1}^{N_x} \eta_j + \zeta
    \end{equation} where $N_x\in\{0,1,\dots\}$ is defined as the number of ``failed" cycles (i.e. those ending at time $\tau(0)$) before a successful one (i.e. one ending at time $\tau(x)$) and has a geometric distribution with success parameter $q(x)\coloneqq \PP_{\nu_x}({\tau}(x)<{\tau}(0))$. For each $x>0$ the laws of the $\eta_j$'s and of $\zeta$ are those of $\tau(0)$ conditioned on the event $\{\tau(0)<\tau(x)\}$ and of $\tau(x)$ conditioned on $\{\tau(x)<\tau(0)\}$, respectively. We remark that although it is not explicit in our notation, the laws of $\eta_1$ and $\zeta$ both depend on $x$.
    
    Lemma \ref{lemma:ruin_proba_order} studies the behavior of $q(x)$ for large $x$; in particular \eqref{eq:q_to_0} shows that it becomes increasingly harder to observe the events of interest as $x$ grows.

    A lower bound for $\PP( \tau^{\nu_x}(x) \leq B(x))$ is given by the probability that $N_x$ in \eqref{eq:sum_decomposition_passage_time} takes an unusually small value. To this end consider for a fixed small $\delta>0$ the events \[E_{\delta}^m \coloneqq \left\{\sum_{j=1}^m \eta_j \leq m(1+\delta)\;\EE\eta_1\right\}\quad \text{and} \quad F_{\delta} \coloneqq \left\{\zeta \leq (1+ \delta)\; \EE\zeta \right\};\] and for each $T>0$: \[n_{\delta}(T) \coloneqq \bigg\lfloor\frac{T-(1+ \delta)\;\EE\zeta}{(1+ \delta)\;\EE\eta_1}\bigg\rfloor.\]
    Then by Lemma \ref{lemma:E_delta_F_delta}: \begin{align*}\label{eq:lower_bound_restart}
      \liminf_{x\to+\infty}  \PP( \tau^{\nu_x}(x) \leq B(x)) & \geq \liminf_{x\to+\infty} \PP\left(\left\{ N_x\leq n_{\delta}(B(x))\right\}\ \cap\ E_{\delta}^{n_{\delta}(B(x))}\ \cap\ F_{\delta}\right)\\
      & = \liminf_{x\to+\infty} \PP\left( N_x\leq n_{\delta}(B(x))\right).
    \end{align*}

    Finally we can conclude the desired uniform lower bound for \eqref{eq:restart_main_result} since:
     \begin{equation*}\label{eq:lower_bpund_restart_II}
        \PP\left( N_x\leq n_{\delta}(B(x))\right) = 1- (1-q(x))^{n_{\delta}(B(x))} = n_{\delta}(B(x))\ q(x)\ (1 + \textbf{o}(1))
    \end{equation*} and $n_{\delta}(B(x))/B(x)$ is bounded away from zero and infinity by Lemma \ref{lemma:asymptotics_tau_0_tau_x}. It follows that for $B(x)$ in the regime considered \[ \PP(\tau(x)\leq B(x))\;B(x)\;q(x)\asymp \PP(\tau(x)\leq B(x))\;B(x)\;\int_0^x\exp(\lambda^* y)\nu_x(dy).\]

The proof of Theorem \ref{thm:restart_main_result_general_measure} is now concluded.
$\bqed$

\end{proof}

\subsection{Proofs of auxiliary lemmas:}
%%% Proof of Lemma 1

\begin{proof}{Proof of Lemma \ref{lemma:ruin_proba_order}:}

We begin with the proof of the upper bound of \eqref{eq:lema_3_statement} by a martingale argument. We have:
\begin{align*}
\int_0^x e^{\lambda^* y} \nu_x(dy) = & \EE_{\nu_x}M_0  =  \EE_{\nu_x}M_{\tau(x)\wedge \tau(0)} \\ = & e^{\lambda^*x}\EE_{\nu_x}[e^{\lambda^*(Z(\tau(x))-x)}\bfone(\tau(x)<\tau(0))]  +\EE_{\nu_x}[e^{\lambda^*Z(\tau(0))}\bfone(\tau(x)>\tau(0))]   \\ \geq & e^{\lambda^*x}\EE_{\nu_x}[\bfone(\tau(x)<\tau(0))] = e^{\lambda^*x}q(x).
\end{align*}
We now turn to the lower bound. Our main ingredient is \cite{Cramer_estimate_Levy_process}\footnote{see \cite{borovkov2013probabilityTheory} for the analogous result for discrete-time random walks}, which states that there exists a constant $C\in (0,1)$ such that \begin{equation}\label{eq:cramer_estimate_LP}
    \lim_{x\to+\infty} e^{\lambda^* x}\PP(\tau(x)<+\infty) = C.
\end{equation} Using that $\tau(0)$ is almost surely finite, for every $x$ and $y\in (0,x)$ \begin{equation}\label{eq:ruin_lower_bound_1}
    \PP_y(\tau(x)<+\infty) = \PP_y(\tau(0)<\tau(x)<+\infty) + \PP_y(\tau(x)<\tau(0)).
\end{equation}

We first claim that \begin{equation}\label{eq:lim_inf_exp_moments_bound}
l = \liminf_{x\to+\infty} \int_0^x e^{\lambda^* y }\nu_x(dy) > 1.  
\end{equation} 

Indeed, $l\geq 1$ and equality cannot hold since otherwise the measures $\nu_x$ would converge (up to taking a subsequence) to the delta measure on $0$, thus contradicting the positive limit of the absorption rates in Lemma \ref{lemma:exponential_rates_converge}, and the claim follows.

In view of the claim above, let $\delta>0$ be small enough such that $C(1-\delta) > C \frac{(1+\delta)}{l}$, and let also $x_{\delta}$ be given by \eqref{eq:cramer_estimate_LP} such that for $x>x_{\delta}$ \[C(1-\delta)\leq e^{\lambda^* x}\PP(\tau(x)<+\infty) \leq C(1+\delta).\]

Noticing that by the Markov property $\PP_y(\tau(0)< \tau(x)<+\infty) \leq \PP_0(\tau(x)<+\infty) = \PP(\tau(x)<+\infty) $, we then have:

\begin{align}
    \begin{split}
        q(x) & = \int_0^x \PP_y(\tau(x)<\tau(0))\nu_x(dy)\\ 
            &=  \int_0^x \PP(\tau(x-y)<+\infty) - \PP_y(\tau(0)< \tau(x)<+\infty) \nu_x(dy) \\
             & \geq \int_0^x \PP(\tau(x-y)<+\infty) - \PP(\tau(x)<+\infty) \nu_x(dy) \\ 
             & \geq  C(1-\delta) e^{-\lambda^*x} \int_0^{x-x_{\delta}} e^{\lambda^* y}\nu_x(dy) + \int_{x-x_{\delta}}^x \PP(\tau(x-y)<+\infty)\nu_x(dy) - C(1+\delta)  e^{-\lambda^*x}\\ 
             & \geq C(1-\delta) e^{-\lambda^*x}\int_0^{x-x_{\delta}} e^{\lambda^* y}\nu_x(dy) +  \PP(\tau(x_{\delta})<+\infty)\nu_x([x-x_{\delta},x]) - C(1+\delta)  e^{-\lambda^*x}.
    \end{split}
\end{align}
It then suffices to find a strictly positive lower bound, independent of $x$ for the following:
\begin{align}\begin{split}\label{eq:lower_bound_3}
    \frac{q(x)e^{\lambda^* x}}{\int_0^x e^{\lambda^* y}\nu_x(dy)} & \geq C(1-\delta)\frac{\int_0^{x-x_{\delta}} e^{\lambda^* y}\nu_x(dy)}{\int_0^x e^{\lambda^* y}\nu_x(dy)} \\ &\hspace{0.9cm} +  \PP(\tau(x_{\delta})<+\infty)\nu_x([x-x_{\delta},x]) \frac{e^{\lambda^* x}}{\int_0^x e^{\lambda^* y}\nu_x(dy)}\\ &\hspace{0.9cm} - C(1+\delta)\frac{1}{\int_0^x e^{\lambda^* y}\nu_x(dy)}.
\end{split}
\end{align}
There are two cases to distinguish: whether $l(\delta) = \liminf\limits_{x\to+\infty} \frac{\int_0^{x-x_{\delta}} e^{\lambda^* y}\nu_x(dy)}{\int_0^x e^{\lambda^* y}\nu_x(dy)}$is equal to $1$ or not. 

In the case $l(\delta)=1$, the first term in the RHS of \eqref{eq:lower_bound_3} dominates the third and the desired lower bound is $C(1-\delta) - C(1+\delta)/l >0$. 
 
If we now assume that $l(\delta) < 1-\varepsilon$ for some positive $\varepsilon$, then the middle term in \eqref{eq:lower_bound_3} is at least $\varepsilon\PP(\tau(x_{\delta})<+\infty)$ for large $x$. The proof concludes by showing that either the first and third terms vanish or the first dominates the third: assume first that there exists a sequence $x_n\to+\infty$ such that \[\int_0^{x_n-x_{\delta}} e^{\lambda^* y}\nu_{x_n}(dy) < 1-\varepsilon' \] for some $\varepsilon'>0$ and all $n$; then $\nu_{x_n}([x_n-x_{\delta},x_n])\geq \varepsilon'$ and consequently \[\int_0^{x_n} e^{\lambda^* y}\nu_{x_n}(dy) \to +\infty. \] For such a sequence, the first and third terms on the RHS of \eqref{eq:lower_bound_3} vanish as $x\to+\infty$. On the other hand, if there is no such sequence (and hence $\liminf \int_0^{x_n-x_{\delta}} e^{\lambda^* y}\nu_{x_n}(dy)\geq 1$), the first term dominates the third one and $\varepsilon\PP(\tau(x_{\delta})<+\infty)$ is still a lower bound for $\eqref{eq:lower_bound_3}$. 

Let us now observe that $q(x)$ goes to $0$ as $x\to+\infty$, so the events of interest are indeed rare events: for any $\varepsilon>0$, there is $k_{\varepsilon}$ such that $\nu[k_{\varepsilon},+\infty) < \varepsilon$; then for $x$ big enough we have that: 
    \begin{align*}
        q(x) &\asymp \exp(-\lambda^*x) \int_0^x \exp(\lambda^*y)\ \nu_x(dy) \\&\leq \exp(-\lambda^*x) \int_0^{k_{\varepsilon}} \exp(\lambda^*y)\ \nu_x(dy)\ +\ \nu_x[k_{\varepsilon},x).
    \end{align*} Since $k_{\varepsilon}$ is fixed the above limit is upper bounded by $\varepsilon$, and $\varepsilon$ is arbitrary. Hence $q(x)\to 0$ as $x\to+\infty$.

The proof of Lemma \ref{lemma:ruin_proba_order} is then complete.
$\wqed$
\end{proof}

%%%%%%%%%%%%%%%%%%%%%%%%%%%%%%%%%%%%%%

\begin{proof}{Proof of Lemma \ref{lemma:asymptotics_tau_0_tau_x}:}
The statement on $\eta$ follows easily from stochastic domination by noticing that  the function $y \mapsto \EE_y \tau(0) $ is non-decreasing. Indeed:  
\begin{align}
\begin{split}
\limsup_{x\to+\infty} \EE[\eta]  = \limsup_{x\to+\infty} \frac{\EE_{\nu_x} [\tau(0)\ \bfone(\tau(0)<\tau(x))]}{1 - q(x)}   \leq \frac{\EE_{\nu} [\tau(0)]}{1 + \textbf{o}(1)}  <+\infty.
\end{split}
\end{align}

%%% \zeta
We now turn to $\zeta$. For each $x>0$ denote by $U(x)$ the overshoot at $x$: $U(x) = Z(\tau(x)) - x$.

by means of the exponential change of measure and the asymptotics of :
 \begin{align*}
 \begin{split}
     \lim_x \EE[\zeta] & = \lim_x \EE_{\nu_x}[\tau(x)\;|\; \tau(x) <  \tau(0)] \\ 
     &=  \lim_x \frac{\EE^{\lambda^*}_{\nu_x}[e^{-\lambda^*U(x)+\lambda^* y_0}\ ;\ \tau(x)\ ;\ \tau(x) <  \tau(0)]}{\EE^{\lambda^*}_{\nu_x}[e^{-\lambda^* U(x) + \lambda^* y_0 }\ ;\ \tau(x) <  \tau(0)]} \\ 
     &\asymp \; \EE^{\lambda^*}_{\nu_{x}}[\tau(x)] \\
     & \asymp \frac{1}{\psi'(\lambda^*)}\int_0^{x} (x - y_0) \; \nu_{x}(dy_0).
 \end{split}
    \end{align*}    

    Using that $\nu$ has finite first moment, the last term above is seen to scale as $x/\psi'(\lambda^*)$ up to multiplicative constants.
\end{proof}\lqqd

%%%%%%%%%%%%%%%%%%%%%%%%%%%%%%%

\begin{proof}{Proof of Lemma \ref{lemma:E_delta_F_delta}:}

    To lighten the notation let us denote $n_{\delta}(B(x))$ by $n_{\delta}$. We will show separately that \begin{equation}\label{eq:lln_eta}
     \lim_x\PP\left(\frac{1}{n_{\delta}}\sum_{j=1}^{n_{\delta}} \eta_j \geq (1+\delta)\; \EE\eta_1 \right) = 0   
    \end{equation}
      and \begin{equation}\label{eq:lln_zeta}
    \lim_x\PP(\zeta > (1+ \delta)\; \EE\zeta) = 0.
    \end{equation} Together they imply \eqref{eq:E_delta_F_delta}.
    
    Let us start with \eqref{eq:lln_eta}. For this part of the proof we apply Chebyshev's inequality:
    
    \begin{align*}
        \lim_x\PP\left(\frac{1}{n_{\delta}}\sum_{j=1}^{n_{\delta}} \eta_j \geq (1+\delta)\; \EE\eta_1 \right) & \leq \lim_x\PP\left(\frac{1}{n_{\delta}^2}\sum_{j=1}^{n_{\delta}} |\eta_j-\EE\eta_j|^2 \geq (\delta\; \EE\eta_1)^2 \right) \\ & \leq \frac{\text{Var}(\eta_1)}{n_{\delta}(\delta\; \EE\eta_1)^2} 
        \\ &  \leq \frac{1}{n_{\delta}(\delta\; \EE\eta_1)^2} \int_0^{+\infty} \EE_{y}[\tau(0)^2]\;\nu(dy)
        \\ &  \asymp \frac{1}{n_{\delta}(\delta\; \EE\eta_1)^2} \int_0^{+\infty} (x-y)^2\;\nu(dy).
    \end{align*} The asymptotic equivalence of the last line above is justified by \cite[Theorem 3]{doney_maller_moments_LP} once we note that $\PP(Z(t)\geq -x_0)$ decays exponentially to $0$ with $t$ for every fixed $x_0>0$ (this is easily checked with an exponential martingale argument using any parameter $\theta>0$ such that $\psi(\theta)<0$). 
    
    Since $\nu$ has finite second moment and does not equal $\delta_0$, $ \EE[\eta_1]$ doesn't go to $0$ with $x\to+\infty$ and we conclude \eqref{eq:lln_eta}.

   Let us now prove \eqref{eq:lln_zeta}. On the one hand the Strong Law of Large Numbers (see \cite[Excercise 7.2 p. 224]{fluctuations_LP_kyprianou}) implies that: \[\lim_{x\to+\infty} \frac{Z(\tau(x))}{\tau(x) \EE^{\lambda^*} Z(1)} =  \lim_{x\to+\infty} \frac{x+U(x)}{\tau(x) \psi'(\lambda^*)} = 1,\; \PP^{\lambda^*}-a.s.\]
    Here as before $U(x)$ denotes the overshoot at $x$. 
    On the other hand, as seen in Lemma \ref{lemma:asymptotics_tau_0_tau_x}, for large $x$, $\EE\zeta$ is asymptotically equivalent to $x/\psi'(\lambda^*)$. Under Cramér's condition $U(x)/x$ converges to $0$ almost surely, so we conclude that \[\lim_{x\to+\infty} \frac{\tau(x)}{\EE^{\lambda^*}[\tau(x)]} = 1,\; \PP^{\lambda^*}-a.s.\] and \eqref{eq:lln_zeta} follows. 
\lqqd

\end{proof}

%%%% Section 

%%%%%%%%%%%%%%%%%%%%%%%%%%%%%%%%%%%%%%%%%%%%%%

\section{Proof of Theorem \ref{thm:main_thm_QSD_Levy}}
\label{sec:proof_restart}

In this section $Z$ denotes either a random walk on discrete time or a L\'evy process, satisfying Cram\'er's condition, and we will consider as restarting measures the family of quasi-stationary distributions $(\nu_x)_{x\geq1}$. The proof of Theorem \ref{thm:main_thm_QSD_Levy} relies on some auxiliary results.

\begin{lemma}\label{lemma:qsd_converge_weakly}
The sequence of QSD measures $(\nu_x)_{x\geq1}$ converges in distribution to the Yaglom limit on $(0,+\infty)$, given by: \[\nu(B) \coloneq \lim_{t\to+\infty} \PP_{x_0}(Z(t)\in B\;|\; t< \tau(0)), \] where $\tau(0) = \inf\{t>0\;:\;Z(t) \leq 0 \}$ and $x_0\in (0,+\infty)$ is an arbitrary point. 
\end{lemma}

Once the sequence of measures admits a weak limit, Remark~\ref{remark:weak_convergence_and_domination}
provides a route to apply Theorem~\ref{thm:restart_main_result_general_measure}, thereby yielding the
lower bound in Theorem~\ref{thm:main_thm_QSD_Levy}. We adopt, however, a different
approach that relies more heavily on properties associated with quasi-stationarity.

% {\color{ree} 

% \begin{remark} If we assume that the QSD for a random walk stochastically dominated by a finite measure $\nu$, then:

% \begin{equation}
%     \beta(x) = -\log \PP_{\nu_x}(\tau(0)\wedge \tau(x)>1) \geq -\log \PP_{\nu_x}(\tau(0)>1).
% \end{equation} Now, since $\PP_y(\tau(0)>1)$ is non-decreasing as a function of $y>0$, $\nu_x\stochdomd \nu $ implies: 
% \begin{equation}
%     \beta(x) \geq -\log \PP_{\nu}(\tau(0)>1).
% \end{equation} Hence $\beta(x)$ has a non-zero limit as $x\to +\infty$ provided that $\PP_{\nu}(\tau(0)>1) \in (0,1)$.   
% \end{remark}
% }

\begin{lemma}\label{lemma:t_x_qsd_is_exponential}
 For the restarted process $Z^{\nu_x}$, $\tau^{\nu_x}(x)$ has an exponential distribution with parameter $\beta(x)q(x)$, where $\beta(x)$ is the exponential rate of absorption associated with $\nu_x$ and $q(x)$ is defined in \eqref{eq:definition_of_q}.
\end{lemma}

\begin{lemma}\label{lemma:exponential_rates_converge}
    The function $\beta(x)$, as defined in Lemma \ref{lemma:t_x_qsd_is_exponential}, is decreasing in $x>0$ and has a positive limit $\beta >0$ as $x\to+\infty$. %In the case of a L\'evy process $\beta$ is the absorption rate associated to the quasi-stationary distribution of the L\'evy process absorbed only at $0$.
\end{lemma}

%%% Proof of THM 3

If we assume the lemmas above we may prove the main result as follows: 

\begin{proof}{Proof of Theorem \ref{thm:main_thm_QSD_Levy}:}\; \\
    By Lemma \ref{lemma:t_x_qsd_is_exponential} we have that $\PP(\tau^{\nu_x}(x)\leq B(x)) = 1 - \exp[-\beta(x)q(x)B(x)]$. Since $B(x)$ grows slower than $e^{\lambda^*x}$ with $x$, $\beta(x)q(x)B(x)\to 0$ and 
        $\PP(\tau^{\nu_x}(x)\leq B(x))\sim\beta(x)q(x)B(x)$ as $x\to+\infty$, with the notation $\sim$ defined at the beginning of Subsection \ref{subsect:Proofs}. 
        
        Also, by Lemma \ref{lemma:exponential_rates_converge}, $\beta(x)$ has a positive limit, and since 
        $\PP(\tau(x)\leq B(x))$ has order $\exp(-\lambda^* x)$, by Proposition \ref{thm:estimates_pasasge_times} we conclude that the ratio \begin{equation}
    \frac{\PP(\tau^{\nu_x}(x)<B(x))}{ \PP(\tau(x)<B(x))}\frac{1}{B(x)\int_0^x e^{\lambda^* y}\nu_x(dy)}
    \end{equation}  is bounded away from zero and infinity uniformly in $x$.
    
    $\bqed$
    
\end{proof}

%%%% Proof of Corollary of THM 3

Before proving the lemmas we compute the quasi-stationary distribution, its exponential moment of order $\lambda^*$ and the absorption rate for a linear Brownian motion explicitly.

\begin{proof}{Proof of Corollary \ref{Coro:linear_BM_QSD}}\label{ex:lin_bm_qsd}

In the case in which $Z$ is distributed as a Brownian motion we can compute all ingredients explicitly. Brownian motion with constant drift $-\mu<0$ has infinitesimal generator: \begin{equation*}
    \cL = \frac{1}{2}\frac{d^2}{dx^2} - \mu\frac{d}{dx}
\end{equation*} whose domain contains the set of twice continuously differentiable functions with right limit (resp left limit) equal to $0$ at $0$ (resp $x$), denoted by $C^2_0(0,x)$. The adjoint operator is: \begin{equation*}
    \cL^* = \frac{1}{2}\frac{d^2}{dx^2} + \mu\frac{d}{dx}
\end{equation*} and has $u(y) = \sin(\pi y/x) e^{-\mu y}$ as a positive eigenfunction. Indeed:  \begin{equation*}
    \cL^* u = -\frac{1}{2}\left(\mu^2 + \frac{\pi^2}{x^2}\right) u.
\end{equation*}

By the spectral characterization of quasi-stationary measures (see e.g. \cite[Proposition 4]{QSD_Villemonais_Meleard}), the measure \begin{equation}\label{eq:qsd_BM}
    \nu_x(dy) = D \sin(\pi y/x) e^{-\mu y} dy,
\end{equation} is then quasi-stationary for $Z$, where $D= \frac{\mu^2x^2+\pi^2}{\pi x (e^{-\mu x} + 1)}$ is the normalization constant. The absorption time $T_x = \tau(x)\wedge \tau(0)$ is exponentially distributed with rate $\beta(x) = \frac{1}{2}\left(\mu^2 + \frac{\pi^2}{x^2}\right)$ that is: 
\begin{equation}\label{eq:qsd_abs_distribution}
    \PP_{ \nu_x}(T_x>t) = e^{-\frac{1}{2}\left(\mu^2 + \frac{\pi^2}{x^2}\right)t}.
\end{equation} 

With the explicit form for the quasi-stationary distribution, and recalling that $\lambda^* = 2\mu$, the exponential moment follows: 
\begin{align}\label{eq:qsd_BM_exp_moment}
    \begin{split}
       \int_0^x e^{\lambda^* y}\nu_x(dy) = D \int_0^x e^{\mu y}\sin\left(\frac{\pi}{x}y\right) dy = D \pi x \frac{e^{\mu x }+ 1}{\mu^2 x^2 + \pi ^2} = e^{\mu x}.
    \end{split}  
\end{align} The result now follows from Theorem \ref{thm:main_thm_QSD_Levy}.
\iffalse

{\color{red}    
    From \eqref{eq:qsd_BM_exp_moment} and Doob’s Optional Stopping Theorem applied to the martingale $M_t = \exp(\lambda^* Z(t))$ we get an explicit expression for $q(x)$:
\begin{equation}
    \begin{aligned}\label{eq:ruin_proba}
       \int_0^x e^{\lambda^* y} \nu_x(dy) =  & \EE_{\nu_x}M_0 = \EE_{\nu_x}M_{\tau(x)\wedge \tau(0)} = e^{\lambda^* x} q(x) + e^{\lambda^* 0}(1-q(x));
    \end{aligned}
\end{equation} so:
    \begin{equation}\label{eq:q(x)_BM}
      q(x) = \frac{\int_0^x e^{\lambda^* y} \nu_x(dy)-1}{e^{\lambda^* x} -1} = \frac{1}{e^{\mu x}+1}.
\end{equation}
}
\fi

\wqed 
\end{proof}

\subsection{Proofs of auxiliary lemmas:}
%%% Proof of Lemma 1
We now turn to the proof of the auxiliary lemmas.

\begin{proof}{Proof of Lemma \ref{lemma:qsd_converge_weakly}:}

Let us recall the piece of notation $T_x = \inf\{t>0\;:\; Z(t) \notin (0,x)\} = \tau(0)\wedge \tau(x)$ and take $x_0$ any number in $(0,+\infty)$. We shall prove that the following limit holds for any bounded measurable function $f$: 
\begin{equation}\label{eq:qsd_convregence_3}
    \lim_{x\to+\infty} \lim_{t\to+\infty} \left|  \EE_{x_0} \left[ f(Z(t)) \; |\; t< \tau(0)  \right] - \EE_{x_0} \left[ f(Z(t)) \; |\; t< T_x  \right]  \right| = 0.
\end{equation} 

The result relies on proving that $\PP_{x_0}(t<\tau(x) \wedge \tau(0)) \sim \PP_{x_0}(t<\tau(0))$ when taking $t$ and $x$ to $+\infty$ in the same order as in \eqref{eq:qsd_convregence_3}. To check this fact, we start by using the strong Markov property: 
\begin{align}\label{eq:qsd_convergence_4}
    \begin{split}
     1 - \frac{\PP_{x_0}(t<T_x) }{\PP_{x_0}(t<\tau(0))} & =  \frac{\PP_{x_0}( \tau(x)\leq t<\tau(0)) }{\PP_{x_0}(t<\tau(0))} \\
      & = \frac{ \EE_{x_0}\left\{\bfone(\tau(x)\leq t)\;\PP_{Z(\tau(x))}\left(t-\tau(x)<\tau(0)\right)\right\}}{\PP_{x_0}(t<\tau(0))}.   
    \end{split}
\end{align}
     
In \cite[Theorem 3.5]{Seva_passage_times_LP_RW} it is proved that for L\'evy processes, for large $t$: \[\PP_{x_0}\left(t<\tau(0)\right) \sim \text{const }V(x_0) \exp(t\;\psi(\lambda_0))\; t^{-3/2},\] where $\lambda_0\in (0,\lambda^*)$ is the unique positive solution to $\psi'(\lambda)=0$, and the function $V$ as defined in \cite[Theorem 2.1]{Seva_passage_times_LP_RW} satisfies: \begin{equation}
     g(x)\; e^{\lambda_0 x} \leq   V(x) \leq e^{\lambda_0 x}\; C
\end{equation} for some constant $C>0$ and some bounded increasing function $g$. It now follows that the right-hand side of~\eqref{eq:qsd_convergence_4} converges to zero. Indeed,
\[
\begin{alignedat}{1}
1 - \frac{\PP_{x_0}(t<T_x)}{\PP_{x_0}(t<\tau(0))}
&\leq
\mathrm{const}\;
\frac{\EE_{x_0}\!\left\{\mathbf{1}(\tau(x)\le t)\,
V(Z(\tau(x)))\right\}}{V(x_0)}
\\[0.8ex]
\overset{t\to+\infty}{\leq}\;&\phantom{\leq\;}
\mathrm{const}\;
\frac{1}{\exp(\lambda_0 x_0)}\,
\EE_{x_0}\!\left\{\mathbf{1}(\tau(x)\le +\infty)\,
\exp\{\lambda_0(x-x_0+U(x))\}\right\}
\\[0.8ex]
\leq\;&\phantom{\leq\;}
\mathrm{const}\;
\EE^{\lambda^*}\!\left\{
\exp\bigl(-(\lambda^*-\lambda_0)(x-x_0+U(x))\bigr)\right\}
\\[0.8ex]
=\;&\phantom{\leq\;}
\mathbf{o}(1),
\hspace{0.3cm} x\to+\infty .
\end{alignedat}
\] Here we used the exponential martingale change of measure and the overshoot $U(x)$ once more (recall Remark \ref{remark:oevrshoot_convergence}). In the case of a random walk under Cram\'er's condition, \cite[Theorem 3.5]{Seva_passage_times_LP_RW} provides an analogous estimate but with a spatial factor $V=1$, so all the above considerations still hold. We may now conclude \eqref{eq:qsd_convregence_3}: 
\begin{align}
\begin{split}
\Bigl|
  \EE_{x_0}\bigl[ f(Z(t)) \mid t<\tau(0) \bigr]
  &-
  \EE_{x_0}\bigl[ f(Z(t)) \mid t<T_x \bigr]
\Bigr|
\\
&=
\frac{
  \Bigl|
    \EE_{x_0}\bigl[ f(Z(t));\, t<\tau(0) \bigr]
    - (1+\mathbf{o}(1))\,\EE_{x_0}\bigl[ f(Z(t));\, t<T_x \bigr]
  \Bigr|
}{
  \PP_{x_0}(t<\tau(0))
}
\\
& \leq \frac{
  \Bigl|
    \EE_{x_0}\bigl[ f(Z(t));\, \tau(x) \le t<\tau(0) \bigr] \Bigr|}{
  \PP_{x_0}(t<\tau(0))
}
    + \mathbf{o}(1)\,\frac{  \Bigl|\EE_{x_0}\bigl[ f(Z(t));\, t<T_x \bigr]
  \Bigr|
}{
  \PP_{x_0}(t<\tau(0))
}
\\
&\le
\|f\|_{\infty}\;
\mathbf{o}(1),
\end{split}
\end{align} and the proof is complete.
\lqqd
\end{proof}

\begin{proof}{Proof of Lemma \ref{lemma:t_x_qsd_is_exponential}:}
We recall that for the process $Z$ absorbed upon exiting $(0,x)$, the absorption time $\tau(x)\wedge\tau(0)$ is exponentially distributed with parameter $\beta(x)$ if started with the quasi-stationary distribution $\nu_x$ (see e.g.\cite[Theorem 2.2]{QSD_book_Collet_etal}). Observe that the $\tau^{\nu_x}(x)$ can be decomposed into a sum of a geometric number of duration of cycles, namely \begin{equation}
    \tau^{\nu_x}(x) = \sum_{j=1}^{N_x} \tau_j,
\end{equation} where $N_x$ is the index of the first cycle at which $Z$ exits through $[x,+\infty)$, and $\tau_j$ is the duration of the $j-$th cycle. Thus $N_x$ is geometric with parameter $q(x) = \PP_{\nu_x}(\tau(x)<\tau(0))$ and the times $\tau_j$ are independent and distributed as $\tau(x)\wedge \tau(0)$. Then $\tau^{\nu_x}(x)$ is a geometric sum of i.i.d. exponential random variables. Moreover, since $\nu_x$ is a quasi-stationary distribution, the time and position at which $Z$ exits $(0,x)$ are independent \cite[Theorem 2.6]{QSD_book_Collet_etal}, and hence so is $N_x$ of $\tau_1,\dots,\tau_{N_x}$.  Being a geometric sum of exponential random variables which are mutually independent and independent of the number of terms, we conclude that $\tau^{\nu_x}(x)$ has an exponential distribution of parameter $\beta(x)q(x)$.
\wqed
\end{proof}

%%% Proof of Lemma 6

\begin{proof}{Proof of Lemma \ref{lemma:exponential_rates_converge}:}
Let $T_x = \tau(x)\wedge \tau(0) = \inf\{t>0:\ Z(t)\notin (0,x)\}$ be the first exit time from $(0,x)$; since the quasi-stationary measure is the unique Yaglom limit for the system under consideration, the rate of absorption $\beta(x)$ is obtained as \begin{equation}
   \beta(x) = -\lim_{t\to+\infty} \frac{1}{t}\log \PP_y(T_x>t)
\end{equation} for any $y\in(0,x)$. Since $x'>x> y >0$ implies $T_{x'}\geq T_x$ $\PP_y -$almost surely, it follows that $\beta(x')\geq \beta(x)$. The sequence of absorption rates is non-increasing with $x$, hence has a limit $\beta$. The same reasoning show that if $\beta_{\infty}$ denotes the absorption rate of the QSD on $(0,+\infty)$ then $\beta(x)\geq \beta_{\infty}$ for every $x\geq 1$ which implies that $\beta\geq\beta_{\infty}>0$. 

%%%% do we put the following?
\iffalse For L\'evy processes we can actually prove that $\beta = \beta_{\infty}$:
To prove that $\beta\coloneqq \lim\limits_{x\to+\infty}\beta(x)$ is strictly positive we will rely on an argument of \cite[Proposition 3.1]{groisman_jonckheere_levy}. For readers' ease let us give a brief review of their result: in \cite{groisman_jonckheere_levy} the authors consider a zero-mean L\'evy process, $\{Z_0(t)\}_{t}$, and characterize the absorption rate of $\{Z_{0}(t)-ct\}_{t}$ at $0$ as $\Gamma(c)\coloneqq\sup_{\theta\in\RR}\{\theta c - \psi_0(\theta)\}$, where $\psi_0$ is the L\'evy exponent of $Z_0$. In our case, $Z$ is a negative-mean L\'evy process, so we may write \[Z(t) = Z_0(t) + t\cdot \EE Z(1),\] where $\EE Z_0(t) = 0$ for all $t$. Next, we observe that $\psi_0(\theta) = \psi(\theta)- \theta\EE Z(1)$ and hence \[\Gamma(c) =  \sup_{\theta\in\RR}\{\theta (c + \EE Z(1)) - \psi(\theta)\} = \zeta[c + \EE Z(1)], \] where $\zeta[\cdot]$ was defined in \eqref{eq:convex_transform}. Finally, we observe that \[\beta = \Gamma(-\EE Z(1)) = \sup_{\theta}\{-\psi(\theta)\}=-\inf_{\theta}\{\psi(\theta)\}.\]
Since $\varphi =\exp\circ \psi $ is strictly convex it has a global minimum attained at $\lambda_0\in (0,\lambda^*)$ and hence $\beta = -\psi(\lambda_0) >0$. Recall that $\lambda_0$ was already mentioned in the paragraph before Section \ref{sec:simulations_restart}.

In the Brownian motion case, we have $\psi(\lambda) =-\mu\lambda + \frac{\lambda^2}{2}$, which gives $\beta  = \frac{1}{2}\mu^2$, as in \eqref{eq:qsd_abs_distribution}.

\fi
$\wqed$
\end{proof}

%%%%%%%%%%%%%%%%%%%%%%%%%%%%%%%%%%%%%%

\section{Conclusion}

This paper investigates various exploration strategies under time constraints in environments with unknown stochastic dynamics focusing on their impact on performance as measured by the time required to reach a set of rare states. 

We aim for this work to be a foundational step towards developing a more qualitative theory of exploration, specifically by incorporating the time needed to observe meaningful signals for the first time. For example in Reinforcement Learning environments with sparse rewards, exploration of the state space or action-state space is widely recognized as a critical bottleneck to efficiency. By addressing this crucial aspect, we seek to contribute to more efficient and effective exploration strategies in such contexts and beyond. 

In this work, we focused on space-invariant one-dimensional dynamics, which provide highly interpretable results and for which we can provide explicit performance guarantees. Building on this understanding, future work will consider more general Markovian dynamics, providing a more comprehensive and realistic framework for exploring challenging environments. \textbf{
A natural extension of this work is to higher-dimensional settings. While explicit hitting-time estimates of the kind used here are generally unavailable in that regime, the ``too-many-particles" phenomenon and the core conclusions of our main results are expected to carry over to considerably more general settings.}

%

%This paper considers different exploration strategies under time constraints for exploring an environment with unknown stochastic dynamics, focusing on their impact on performance measured by the time required to reach a set of rare states. In Reinforcement Learning, especially in environments with sparse rewards, exploration — whether of the state space or the action-state space — is well identified by the community as a significant bottleneck to efficiency. We expect this work to provide a first step toward developing a more qualitative theory of exploration, specifically by incorporating the time needed to observe meaningful signals for the first time. By addressing this critical aspect, we hope to contribute to more efficient and effective exploration strategies, ultimately advancing the understanding and application of Reinforcement Learning in challenging environments with rare events.
%In future work, more generic Markovian dynamics will be considered.

%\THEEndNotes
%\begingroup \parindent 0pt \parskip 0.0ex \def\enotesize{\normalsize} \theendnotes \endgroup

% Appendix here
% Options are (1) APPENDIX (with or without general title) or
%             (2) APPENDICES (if it has more than one unrelated sections)
% Outcomment the appropriate case if necessary
%
% \begin{APPENDIX}{<Title of the Appendix>}
% \end{APPENDIX}
%
%   or
%
% \begin{APPENDICES}
% \section{<Title of Section A>}
% \section{<Title of Section B>}
% etc
% \end{APPENDICES}

% Acknowledgments here
\ACKNOWLEDGMENT{P.B. and E.G. would like to express their deep gratitude to professors E. Mordecki and J.R. Le\'on for very valuable discussions.}

% References here (outcomment the appropriate case)

% CASE 1: BiBTeX used to constantly update the references
%   (while the paper is being written).
\bibliographystyle{informs2014} % outcomment this and next line in Case 1
\bibliography{biblio} % if more than one, comma separated

%\bibliographystyle{informs2014} % outcomment this and next line in Case 1
%\bibliography{sample} % if more than one, comma separated

% CASE 2: BiBTeX used to generate mypaper.bbl (to be further fine tuned)
%\input{mypaper.bbl} % outcomment this line in Case 2

%If you don't use BiBTex, you can manually itemize references as shown below.

%\bibliographystyle{nonumber}

\section*{Appendix I: M/M/1}
\label{appendix:MM1}

In this section we investigate the interplay between the probability of exploration, under the time
regime of interest, and the efficiency of the associated estimation procedure.  
In particular, we show that the variance of the natural Monte Carlo estimator mirrors the typical exploration time.

Consider an $\mathrm{M/M/1}$ queue on the truncated state space $\{0,1,\dots,n\}$ with load parameter $\rho<1$, and assume that $X(0)=0$ almost surely.  
The chain is ergodic, and its invariant measure is geometrically decaying: for any $k\in\mathbb{N}$ with $k\le n$, the stationary weight is proportional to $\rho^k$.  
By the Ergodic Theorem, the stationary measure~$\pi$ satisfies
\[
\pi(k)
    = \lim_{t\to\infty} \frac{1}{t}\int_0^t \mathbf{1}\bigl(X(s)=k\bigr)\,ds .
\]
A Monte Carlo estimator for $\pi(k)$, based on a time horizon $B(k)$, is therefore
\[
    \widehat{\pi}(k)
        = \frac{1}{B(k)}\int_0^{B(k)} \mathbf{1}\bigl(X(s)=k\bigr)\,ds .
\]

We focus on restrictive time budgets and assume that the estimation of $\pi(k)$ for large $k$ is performed with a linear-in-$k$ horizon,
\[
    B(k)\asymp k, \qquad k\to\infty .
\]
In such a regime, the process is expected to spend exceedingly little time in state~$k$: indeed, the expected hitting time of~$k$ grows exponentially fast in~$k$.  
As we show next, this scarcity of observations is the principal obstruction to reliable estimation, and it directly governs the variance of the Monte Carlo estimator.

Fix $c>0$ and introduce the surrogate estimator
\[
    \xi(k) \coloneqq c\, \mathbf{1}\bigl(\tau(k)< B(k)\bigr),
\]
where $\tau(k)$ denotes the hitting time of state~$k$. \textbf{Then, by noting that $\widehat{\pi}(k) + \xi(k)$ is almost surely bounded by $1+c$, and that $\bigl| \widehat{\pi}(k) - \xi(k) \bigr|\bfone(\tau(k)>B(k)) = 0$, one obtains: }

\begin{multline*}
\bigl|\operatorname{Var}\,\widehat{\pi}(k) - \operatorname{Var}\,\xi(k)\bigr|
  = \Bigl|\mathbb{E}_0\bigl[(\widehat{\pi}(k)-\xi(k))(\widehat{\pi}(k)+\xi(k))\bigr] \\
   - \mathbb{E}_0(\widehat{\pi}(k)-\xi(k))\,
    \mathbb{E}_0(\widehat{\pi}(k)+\xi(k))\Bigr| \\
  \leq 2(1+c)\cdot \mathbb{E}_0|\widehat{\pi}(k) -\xi(k)| \\
  \leq \mathrm{const}\cdot \mathbb{P}_0\!\bigl(\tau(k)\le B(k)\bigr),
\end{multline*} for some constant independent of $k$, thus showing that the variance is essentially governed by the rare-event probability of reaching~$k$ within the available time budget.

\section*{Appendix II: The Fleming-Viot particle system and FVRL}
\label{appendix:FVRL}

The approach found in \cite{FV_Jonckheere_etal} for defining a restarted mechanism is based on the paradigm of Fleming-Viot (FV) particle system, introduced by \cite{Burdzy_FV_original}. In this framework, a population of particles (which can represent different simulations in parallel) evolves over time, exploring the state space in parallel. When a particle falls into a less promising region $ A$, it is ``restarted" by being replaced with a copy of another, more promising particle. This ensures that resources are concentrated on exploring the most fruitful areas of the state space. The FV strategy is particularly effective in scenarios where certain regions of the state space are more likely to yield valuable rewards, as it dynamically reallocates exploration effort towards these regions, thereby increasing the overall efficiency of the exploration process.

Note that various papers \cite{Asselah_Ferrari_Groisman_QSd_finite_spaces, ferrari_maric, Villemonais_FV_to_QSD_diffusions} have shown that FV
empirical measures converge as the number of particles tend to infinity to the conditioned evolution of the process, i.e.,
$$ m_t^N(A) \to P(X_t \in A | T \ge t ), N \to \infty,$$
where $T$ is the hitting time of $ A$. On the other hand $P(X_t \in A | T \ge t ) \to \nu(A)$, as $t \to \infty$ where $\nu$ is the (a for countable state space) QSD associated with the process absorbed in $ A$. This is the motivation for restarting according to the  QSD in our work.

\section*{Appendix III: An example in RL: the control of blocking in a $M/M/1/K$}
\label{appendix:completeRL}
This example has been considered in \cite{FV_Jonckheere_etal} and generalized to multidimensional settings (several coupled queues).
It serves here as a canonical illustration of Reinforcement Learning in environments with \emph{sparse and rare rewards} where the exploration is the main bottleneck. We summarize the model and the FVRL strategy introduced in \cite{FV_Jonckheere_etal} as a motivation for our theoretical results.
Consider an $M/M/1/K$ queue with \textbf{fixed} buffer capacity $K$, arrival rate $\lambda$, service rate $\mu$, and load $\rho = \lambda/\mu < 1$. The state space is $\cS = \{0,1,\dots,K\}$, representing queue occupancy. The control action $a \in \{0,1\}$ determines whether to accept ($a=1$) or block ($a=0$) arriving jobs. The reward function is structured to penalize both the under-utilisation of \textbf{resource} and large blocking probabilities:
$$
r(x,a) =
\begin{cases}
0 & \text{if } a = 1 \\
B(1 + b^{x - x_{\text{ref}}}) & \text{if } a = 0
\end{cases}
$$
with $B > 0$, $b > 1$, and reference state $x_{\text{ref}}$.
This structure induces threshold-optimal policies (blocking all the incoming \textbf{traffic} after a given threshold).
\textbf{For threshold policies parameterized by $\theta$, we define the policy-dependent blocking threshold $K_\theta = \lfloor \theta \rfloor + 1$, representing the state at which the policy $\pi_\theta$ begins rejecting incoming jobs. To ensure that the blocking threshold never exceeds the physical buffer, the policy parameter is restricted to the compact domain $\Theta = [0, K-1]$, so that $K_\theta \leq K$ is guaranteed at all times.}
Under the average reward criterion and using a policy gradient strategy, one aims at minimizing:
$$
J^{\pi_\theta} = \sum_{x \in \cS} p^{\pi_\theta}(x) \sum_{a} \pi_\theta(a|x) Q_{\theta}(x,a),
$$
where $\theta \in \Theta$ is the parametrisation of the policy. The policy gradient theorem yields:
\begin{equation}
\nabla_\theta J^{\pi_\theta} = \sum_{x \in \cS} p^{\pi_\theta}(x) \sum_{a} Q_{\theta}(x,a) \nabla_\theta \pi_\theta(a|x)
\end{equation}
For threshold policies parameterized by $\theta \in \Theta$, this simplifies to:
\begin{equation} \label{grad:est}
\nabla_\theta J^{\pi_\theta} \propto p^{\pi_\theta}(\mathbf{K_\theta}-1) \left[ Q_{\theta}(\mathbf{K_\theta}-1,1) - Q_{\theta}(\mathbf{K_\theta}-1,0) \right]
\end{equation}
\textbf{where $K_\theta = \lfloor \theta \rfloor + 1$ is the learned blocking threshold, and the gradient is concentrated at state $K_\theta - 1 = \lfloor \theta \rfloor$, the last state before the policy blocks under the current parameterization $\theta \in \Theta$.}

\subsubsection*{The Gradient Estimation Problem:}

Given this framework, the RL strategy will be efficient if the gradient 
estimates can be informative. We discuss here only the estimate of the 
stationary probability $p^{\pi_\theta}(\mathbf{K_\theta}-1)$.
The stationary distribution for state $\mathbf{K_\theta}$, under the 
threshold policy $\pi_\theta$\textbf{, is given by the $M/M/1/K_\theta$ 
formula:}
$$
p^{\pi_\theta}(\mathbf{K_\theta}) = 
\frac{(1-\rho)\rho^{\mathbf{K_\theta}}}{1-\rho^{\mathbf{K_\theta}+1}}.
$$
For typical parameters ($\rho = 0.7$, \textbf{$K_\theta = 40$}), 
$p^{\pi_\theta}(\mathbf{K_\theta}) \sim \mathcal{O}(10^{-7})$. Since 
$p^{\pi_\theta}(\mathbf{K_\theta}-1) \sim 
p^{\pi_\theta}(\mathbf{K_\theta})/\rho$, both are exponentially small.

Consequently:
\begin{itemize}
\item Vanilla Monte Carlo estimators of 
$p^{\pi_\theta}(\mathbf{K_\theta}-1)$ exhibit prohibitive variance 
(in terms of multiplicative errors),
\item Much more importantly, the policy gradient in~\eqref{grad:est} 
becomes numerically zero under finite sampling. Hence, the learning 
signal vanishes, preventing convergence to an optimal $\theta$.
\end{itemize}

\subsection{Parallel sampling}
\begin{figure}
    \centering
    \includegraphics[width=\linewidth]{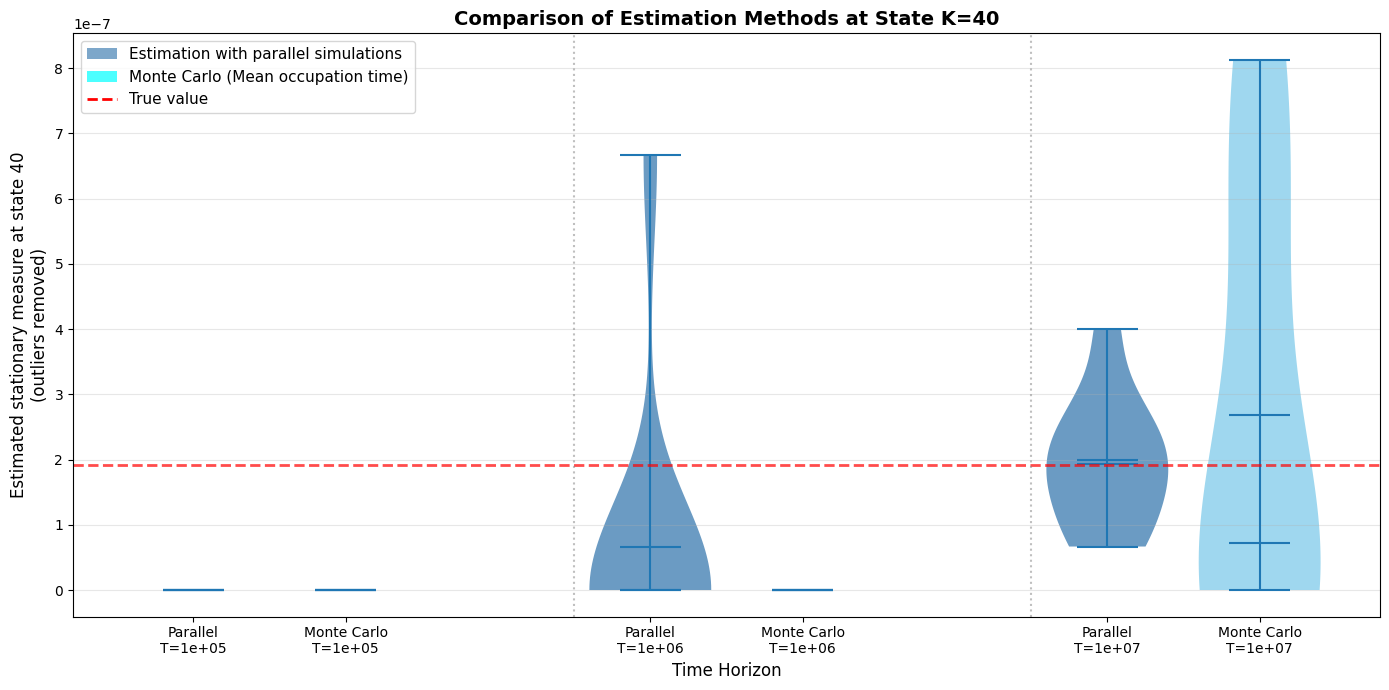}
    \caption{Simulation of parallel $\mathrm{M/M/1}$ queues with $K=40$, $J=12$, $\lambda=0.7$, and $\mu=1$, comparing renewal-based stationary probability estimates across time horizons $10^5$, $10^6$, and $10^7$. The number of parallel copies to deploy in each time regime is computed using the results of our main theorems.}
    \label{fig:parallel_MM1}
\end{figure}

In order to compare our results with those of \cite{FV_Jonckheere_etal}, 
we present a set of simulations for parallel $\mathrm{M/M/1}$ queues under 
different time regimes.  
The hyperparameters used in Figure~\ref{fig:parallel_MM1} are: total queue 
capacity $K=40$, initial condition $J=12$, arrival rate $\lambda=0.7$, and 
service rate $\mu=1$. \textbf{In this numerical illustration, $K_\theta = K 
= 40$, so that $p^{\pi_\theta}(K)$ is the exponentially rare 
quantity identified in the Gradient Estimation subsection above.}

Following \cite{FV_Jonckheere_etal}, we estimate the stationary probability 
via a renewal-type decomposition based on successive returns to $J$.  
More precisely, we estimate (i) the expected return time to state~$J$, 
which in our setting is of order \(10^2\), and independently (ii) the 
probability of hitting~$K$ \emph{starting from~$J$}.  
For the latter we report estimators under time horizons \(B(K)\) of 
magnitudes \(10^5\), \(10^6\), and \(10^7\), and we compare them with the 
true stationary value \(\mathbf{p^{\pi_\theta}(K)}\) shown in 
Figure~\ref{fig:parallel_MM1}.

To ensure a fair comparison with the experiment of 
\cite[Section~4.1]{FV_Jonckheere_etal}, we remark that the total number of 
events observed in our simulations---counting both arrivals and 
departures---is one to two orders of magnitude larger than theirs.  
Indeed, across our time regimes we record approximately \(1.7\times 10^r\) 
events for \(r=5,6,7\), which exceeds the event counts reported in 
\cite[Figure~2]{FV_Jonckheere_etal}.  
This discrepancy is natural: correlations in the Fleming--Viot particle 
system substantially facilitate the exploration of the rare state~$K$.

In summary, and as expected, the estimation accuracy with independent 
particles lies between that of naive Monte Carlo (which rarely reaches $K$ 
even for time of order $10^6$) and that of the Fleming--Viot particle system.

\subsection{Fleming-Viot Solution}

The FVRL method addresses this by introducing an \emph{absorption set}, where prior knowledge is used on the fact that no rewards is granted in $ \A = \{0,1,\dots,J-1\} $ with $ J < K $. The Fleming-Viot particle system consists of $ N $ particles $ (\xi_t^{\nu}(i))_{i=1}^N $ evolving in $ \A^c $, with resetting mechanism: particles hitting $ \A $ jump to positions of randomly selected surviving particles. This particle system approximately estimates the quasi-stationary distribution:
\[
\nu_Q^{\pi}(x) = \lim_{t\to\infty} \PP^{\partial\A^c}(X_t^{\pi} = x | T_{\mathcal{K}} > t)
\]
where $ T_{\mathcal{K}} $ is the hitting time of $ \A $.
The FV-RL algorithm leverages this to construct gradient estimators:
\begin{equation}
\widehat{\nabla_\theta v^{\pi_\theta}} = \sum_{x \in \AA^c} \hat{p}^{\pi_\theta}(x) \sum_{a} \hat{Q}_{\theta}(x,a) \nabla_\theta \pi_\theta(a|x)
\end{equation}
where $ \hat{p}^{\pi_\theta} $ is estimated via the FV particle system.

For $ K = 40 $, $ \rho = 0.7 $:
\begin{itemize}
\item Vanilla MC estimates $ p^{\pi}(K) \approx 0 $ even with $ 10^6 $ samples
\item FV provides accurate estimates with $ N \sim 10^3 $ particles
\item FVRL converges to $ K^* $ while MC-based policy gradient fails completely
\end{itemize}

The method effectively trades the rare event probability $ p^{\pi}(K) $ for the substantially larger quasi-stationary probability $ \nu_Q^{\pi}(K) $, overcoming the exponential sample complexity of vanilla approaches. We do not compare
here with our results of the estimates 
constructed in \cite{FV_Jonckheere_etal} are more complex than the idealized situation described in Section 3.4 which can however serve as a rule of thumb for practitioners.

\end{document}